\documentstyle[amssymb,12pt]{amsart}
\newtheorem{prop}{Proposition}[section]
\newtheorem{defin}{Definition}[section]
\newtheorem{lemma}{Lemma}[section]
\newtheorem{theorem}{Theorem}[section]
\newtheorem{corollary}{Corollary}[section]

\theoremstyle{remark}
\newtheorem{remark}{Remark}
\newtheorem{notation}{Notation}

\newtheorem{convention}{Convention}

\begin{document}
\newcommand{\nc}{\newcommand}
\nc{\pa}{\partial}
\nc{\cA}{{\cal A}}\nc{\cB}{{\cal B}}\nc{\cC}{{\cal C}}
\nc{\cE}{{\cal E}}\nc{\cG}{{\cal G}}\nc{\cH}{{\cal H}}
\nc{\cX}{{\cal X}}\nc{\cR}{{\cal R}}
\nc{\Ad}{\on{Ad}}\nc{\sh}{\on{sh}}
\nc{\ad}{\on{ad}}\nc{\Der}{\on{Der}}\nc{\End}{\on{End}}\nc{\res}{\on{res}}
\nc{\Imm}{\on{Im}}\nc{\limm}{\on{lim}}
\nc{\de}{\delta}\nc{\si}{\sigma}\nc{\ve}{\varepsilon}
\nc{\on}{\operatorname}
\nc{\al}{\alpha}
\nc{\CC}{{\mathbb C}}\nc{\ZZ}{{\mathbb Z}}\nc{\NN}{{\mathbb N}}
\nc{\AAA}{{\mathbb A}}\nc{\cO}{{\cal O}}\nc{\cK}{{\cal K}}
\nc{\la}{{\lambda}}\nc{\G}{{\mathfrak g}}\nc{\A}{{\mathfrak a}}
\nc{\SS}{{\mathfrak S}}
\nc{\HH}{{\mathfrak h}}
\nc{\N}{{\mathfrak n}}
\nc{\La}{\Lambda}
\nc{\g}{\gamma}\nc{\eps}{\epsilon}\nc{\wt}{\widetilde}
\nc{\wh}{\widehat}

%
%
%

\newcommand{\ldar}[1]{\begin{picture}(10,50)(-5,-25)
\put(0,25){\vector(0,-1){50}}
\put(5,0){\mbox{$#1$}} 
\end{picture}}

\newcommand{\lrar}[1]{\begin{picture}(50,10)(-25,-5)
\put(-25,0){\vector(1,0){50}}
\put(0,5){\makebox(0,0)[b]{\mbox{$#1$}}}
\end{picture}}

\newcommand{\luar}[1]{\begin{picture}(10,50)(-5,-25)
\put(0,-25){\vector(0,1){50}}
\put(5,0){\mbox{$#1$}}
\end{picture}}

\title[Quasi-Hopf algebras associated with ${\mathfrak{sl}}_{2}$ and complex
curves]{Quasi-Hopf algebras associated with ${\mathfrak{sl}}_{2}$ and complex
curves}

\author{Benjamin Enriquez}
\address{Centre de Math\'ematiques, URA 169 du CNRS, Ecole Polytechnique,
91128 Palaiseau, France}

\address{FIM, ETH-Zentrum, HG G44, R\"amistr. 101, CH-8092 Z\"urich,
  Switzerland}

\author{Vladimir Rubtsov}
\address{V.R.: D\'ept. de Math\'ematiques, Univ. d'Angers, 
2, Bd. Lavoisier, 49045 Angers, France}

\address{ITEP, 25, Bol. Cheremushkinskaya, 117259 Moscou, Russia}

\date{August 1996; revised May 1998}

\begin{abstract} We construct quasi-Hopf algebras quantizing double
extensions of the Manin pairs of Drinfeld, associated to a curve with a
meromorphic differential, and the Lie algebra ${\mathfrak{sl}}_{2}$. 
This construction makes use of an analysis of the vertex relations for
the quantum groups obtained in our earlier work, PBW-type results and
computation of $R$-matrices for them; its key step is
a factorization of the twist operator relating ``conjugated'' versions
of these quantum groups. 
\end{abstract}

\maketitle

\subsection*{Introduction}
In \cite{D}, V. Drinfeld introduced Manin pairs, attached to an
absolutely simple
Lie algebra over a complex curve $X$ with a meromorphic
differential $\omega$. In the case where the genus of
$X$ is zero or one, special cases of his construction give rise to Manin
triples, whose quantization are the Yangians, quantum affine algebras and
elliptic algebras. Drinfeld posed the problem of quantizing these
general Manin pairs, in the sense of quasi-Hopf algebras. 

In this paper, we solve this problem in the untwisted case where the Lie
algebra is equal to $\A\otimes \CC(X)$, 
$\A={\mathfrak{sl}}_{2}$, for an arbitrary curve $X$.  

In that case, Drinfeld's Manin pair presents itself as follows. Let $S$
be a finite 
set of points of $X$ containing the zeroes and poles of $\omega$,
$k_{s}$ be the local field at each $s\in S$, and $k=\oplus_{s\in
S}k_{s}$. Endow $\A \otimes k$ with the scalar product given by the
tensor product of the Killing form of $\A$ and $\langle f,g \rangle_{k}
= \sum_{s\in S}\res_{s}(fg\omega)$. Let $R$ be the ring of functions on
$X$, regular outside $S$; it can be viewed as a subring of $k$. 
The ring $R$ is a Lagrangian (that is, maximal isotropic) subspace of
$k$ (a proof of this fact is in the Appendix), so that 
$(\A\otimes k, \A\otimes R)$ forms a Manin pair. 

In our earlier paper \cite{ER}, we introduced a double extension
$(\G,\G_{R})$ of this Manin pair. The Lie algebra $\G$ is a direct sum
$(\A\otimes k) \oplus \CC D \oplus \CC K$, with $K$ a central element and $D$
a derivation element, and $\G_{R}$ is equal to $(\A\otimes R)\oplus \CC
D$. In \cite{ER}, we considered a certain Manin triple
$(\G,\G_{+},\G_{-})$, obtained from this pair by a classical twist, and
we constructed a quantization $U_{\hbar}\G$ of this Manin triple. 

Let us describe in more detail the Manin triple of \cite{ER}. 
Let $\A=\N_{+}\oplus \HH \oplus \N_{-}$ be a Cartan decomposition of
$\A$. Let $\Lambda$ be a Lagrangian complement of $R$ in $k$. 
The Lie algebra
$\G_{+}$ is then defined as $(\HH\otimes R) \oplus (\N_{+}\otimes k) \oplus
\CC D$, and $\G_{-}$ as $(\HH\otimes \La) \oplus (\N_{-}\otimes k) \oplus
\CC K$. 

The Manin triple $(\G,\G_{+},\G_{-})$ has the following interpretation. 
Recall that the affine Weyl group of $\G$ is the semi-direct product of
a group of translations, isomorphic to $\ZZ^{S}$, with the Weyl group of
$\A^{S}$. The triple $(\G,\G_{+},\G_{-})$ can
then be viewed as the limit, for the
length of the Weyl group element becoming infinite, of the triple
obtained from $(\G,\G_{R}, \G_{\La})$ [with $\G_{\La}= 
(\A\otimes \La)\oplus \CC
K$] by conjugation by an affine Weyl group element corresponding to a
positive translation. 
A similar procedure had been employed earlier by Drinfeld in
\cite{new-real}. We recall the results of \cite{ER} in sections
\ref{mpt} and \ref{qmt}. 

Let us now describe the main points of the present work. We first
construct a subalgebra $U_{\hbar}\G_{R} \subset U_{\hbar}\G$; this
inclusion deforms the inclusion of Lie algebras $\G_{R}\subset \G$
(section \ref{UgR}).  This construction is as follows.  The main
difficulty of \cite{ER} was to produce the correct relations for the
quantum counterparts of the algebras generated by $\N_{\pm}\otimes k$.
These relations are presented in terms of generating series; they are
usually called vertex relations.  We observe that there exists a
system of such relations, in which the products of generating series
are multiplied by scalar functions belonging to $R\otimes R$. This
system can be modified so as to define vertex relations for an algebra
$U_{\hbar}\G_{R}$ (sect. \ref{UgR.1}), which turns out to be a
deformation of the enveloping algebra $U\G_{R}$ (section \ref{UgR.2}).
To show that this algebra is indeed a subalgebra of $U_{\hbar}\G$, we
have to establish Poincar\'e-Birkhoff-Witt (PBW) type results for
$U_{\hbar}\G$ (Prop. \ref{PBW-for-Ug}). These results follow from
general similar results on algebras presented by vertex relations
(section \ref{vx}) and rely isomorphism with a formal version of the
Feigin-Odesski shuffle algebras (see \cite{FO}).

Let $\Delta$ denote the coproduct of $U_{\hbar}\G$. The subalgebra 
$U_{\hbar}\G_{R}$
then satisfies $\Delta(U_{\hbar}\G_{R})$ $\subset 
U_{\hbar}\G\hat\otimes U_{\hbar}\G_{R}$ (see Prop. \ref{additional}).
This motivates the following construction.

Consider the Manin triple $(\G,\bar{\G}_{+},\bar{\G}_{-})$, obtained
from $(\G,\G_{R}, \G_{\La})$ as the limit of conjugations by negative
translations (or equivalently, from $(\G,\G_{+},\G_{-})$ by
conjugation by $\on{card}(S)$ copies of the nontrivial element of the
Weyl group of $\A$).  Using the results of \cite{ER}, it is easy to
produce a quantization $(U_{\hbar}\bar{\G},\bar\Delta)$ of this Manin
triple. We then show that the Hopf algebras $U_{\hbar}\G$ and
$U_{\hbar}\bar{\G}$ are isomorphic as algebras (Prop. \ref{isom}), and
their coproducts are conjugated under some $F\in
U_{\hbar}\G^{\hat\otimes 2}$ (Prop. \ref{prop:F:twist}), which also
satisfies cocycle identities (Prop. \ref{prop:F:cocycle}); in other
words, both Hopf algebras are connected by a twist, in the sense of
\cite{D}.  Let us explain how this result is obtained.

$F$ is constructed as $\sum_i \al^i \otimes \al_i$, for
$(\al^i),(\al_i)$ dual bases of the subalgebras $U_{\hbar}\N_+$ and
$U_{\hbar}\N_-$ of $U_{\hbar}\G$ generated by the deformations of
$\N_{\pm} \otimes k$. The universal $R$-matrices of
$(U_{\hbar}\G,\Delta)$ and $(U_{\hbar}\bar{\G},\bar\Delta)$ are then
expressed simply as the products of $F$ and a factor $\cK$ depending
on the Cartan modes (Prop.  \ref{prop:R:mat}). One checks that $\cK$
is a twist connecting $\bar\Delta$ and $\Delta'$ (Lemmas
\ref{K:cocycle} and \ref{K:conj}). ($\Delta'$ is $\Delta$ composed with
the permutation of factors.)  Therefore, since $\cR^{-1}$ is a twist
connecting $\Delta'$ and $\Delta$, it follows that $F$ is a twist
connecting $\Delta$ and $\bar\Delta$ (Props. \ref{prop:F:cocycle} and
\ref{prop:F:twist}).

On the other hand, we have $\bar\Delta(U_{\hbar}\G_{R})\subset 
U_{\hbar}\G_{R}\hat\otimes U_{\hbar}\G$ (Prop. \ref{additional}).

It is then easy to see that any factorization of $F$ of the form
$F=F_{2}F_{1}$, $F_{1}\in U_{\hbar}\G \otimes U_{\hbar}\G_{R}, 
F_{2}\in U_{\hbar}\G_{R} \otimes U_{\hbar}\G$,
yields an algebra morphism $\Delta_{R}$
from $U_{\hbar}\G_{R}$ to
$U_{\hbar}\G_{R}^{\hat\otimes 2}$, by the formula $\Delta_{R}=\Ad(F_{1})
\circ \Delta$. 
Moreover, the associator $\Phi$ of 
the quasi-Hopf algebra obtained from $U_{\hbar}\G$ by the twist by
$F_{1}$, belongs to $U_{\hbar}\G_{R}^{\hat\otimes 3}$
(Prop. \ref{Phi-in-R}). A simple argument shows that the antipode
$S_{R}$ of $U_{\hbar}\G$ corresponding to this twist, 
preserves $U_{\hbar}\G_{R}$. This shows that $(U_{\hbar}\G_{R},
\Delta_{R},\Phi, S_{R})$ is a
sub-quasi-Hopf algebra of the twist by $F_{1}$ of $(U_{\hbar}\G,\Delta)$. 

To obtain a factorization of $F$ is thus the key point of our
construction. This is achieved in section \ref{fctr}. Any possible 
solution $(F_{1},F_{2})$ of the factorization identity 
is expressed simply in terms of 
left and right $U_{\hbar}\G_{R}$-module maps $\Pi,\Pi'$ from 
$U_{\hbar}\G$ to $U_{\hbar}\G_{R}$ applied to $F$, and of variable
elements of $U_{\hbar}\G_{R}^{\otimes 2}$. The difficulty is to show that for
some choice of those elements, the factorization identity is
satisfied. This is equivalent to showing that $F^{-1}[(\Pi \otimes 1)F]$
and $[(1\otimes \Pi')F]F^{-1}$ belong to $U_{\hbar}\G\otimes
U_{\hbar}\G_{R}$, resp. to $U_{\hbar}\G_{R}\otimes
U_{\hbar}\G$. For this, we use the pairing between
the quantizations $U_{\hbar}\G_{+}$ and $U_{\hbar}\G_{-}$ of $\G_{+}$
and $\G_{-}$, and the computation of the orthogonals of their intersections
with $U_{\hbar}\G_{R}$ (Prop. \ref{orthBR}).  

We close the paper by some remarks related to our
construction. We observe that the quasi-Hopf algebras $U_{\hbar}\G$ and
$U_{\hbar}\G_{R}$  
fit in an inductive system w.r.t. the relation $S\subset S'$, 
and that the corresponding inductive limit is a quantization of double
extensions of the adelic versions of Drinfeld's Manin pairs (section
\ref{adal}). We also find an algebra automorphism of $U_{\hbar}\G$,
deforming the action of the generator of the Weyl group of $\A$ (section
\ref{qwga}). Section \ref{agen} is devoted to analogues and
generalizations of $U_{\hbar}\G$. In \ref{galaut}, we exploit the fact
that the central terms occur only in the exponential form
$\exp(\hbar K \partial)$ [where $\partial$ is the derivation of $k$
defined by $\partial f = df / \omega$]  
to construct analogues of $U_{\hbar}\G$, where $\exp(\hbar K \partial)$
is replaced by a more general automorphism of $k$.  In \ref{disc}, we
construct analogues of those algebras and of their Weyl group
automorphism, associated with discrete sets. 

An expression of $F$ was given in an earlier version of this work.
However, this expression is not correct, as it was pointed out in
\cite{DK}. In Remark \ref{correction}, we discuss this problem and
how this modifies the proofs of the results of \cite{EF1,EF2}, which
remain valid.

Let us now mention some possible extensions of this work. It is natural
to ask how the algebras introduced here depend on the pair $(X,
\omega)$. The algebra $U_{\hbar}\G$ probably possesses level $1$
modules similar to those studied in the Yangian and quantum affine
cases. It would then be interesting to study quantum
Knizhnik-Zamolodchikov type equations for traces of corresponding
intertwining operators. 
Another subject of interest could be the representation theory of
$U_{\hbar}\G_{R}$. In \cite{ER}, we studied level zero representations
of $U_{\hbar}\G$, indexed by formal discs around each point of $S$;
these representations are also $U_{\hbar}\G_{R}$-modules, and as such
their parameter could probably take values outside those discs. 

Finally, the question arises whether the formulas defining
$U_{\hbar}\G$ can be written is closed form (rather than in the sense
of formal series) and can be analytically continued to complex values
of $\hbar$.  In general, the solution to this problem might be related
to functional equations satisfied by the structure constants of this
algebra. Let us however mention two cases where the answer to this
question is positive. One of them is when $X$ is an elliptic curve,
and $\omega = dz$. This case was treated by G. Felder and one of us
(\cite{EF1}). We showed, using arguments of the present paper, the
connection of the algebra $U_{\hbar}\G$ with the elliptic quantum
groups of \cite{F}. The other case was treated in \cite{Some}. There
we study the case of a genus $>1$ curve $X$, with differential
$\omega$ regular and having only double poles. In that case, the
structure constants of $U_{\hbar}\G$ involve some theta-functions and
odd theta-characteristics of $X$.  Also let us mention the work
\cite{DI}, where ``analytic'' algebras with close analogy to
$U_{\hbar}\G$ were introduced.

We would like to express our thanks to 
V. Fock, C. Fronsdal, S. Khoroshkin, D. Lebedev, S. Majid,
N. Reshetikhin and M. Semenov-Tian-Shansky for discussions about this
work. The second author is supported by grant RFFI-96-01-01101.

\medskip\section{Manin pairs and triples} \label{mpt}
\
This shows that subsection{Completed tensor products and algebras}

Let $V,W$ be complex Tate's vector spaces (see \cite{BFM} for the
definition; the
only examples of Tate's vector spaces we will use are either discrete or
isomorphic to the
sum of a finite number of copies of a field of Laurent power series in
one variable). 

The completed tensor product of $V$ and $W$ is defined as the inverse
limit 
$$
\lim_{\leftarrow a,b} 
(V \otimes W/V_{a} \otimes W_{b}),
$$ 
$V_{a},W_{b}$ being a system
of vector spaces, that consitute neighborhoods of zero in $V,W$, and
denoted by $V\hat\otimes W$. For example, with $V = \CC((v))$ and $W =
\CC((w))$ endowed with the $v$- and $w$-adic topologies, we have $V
\hat\otimes W = \CC[[v,w]][v^{-1},w^{-1}]$ which we will denote as
$\CC((v,w))$. 

With the same notation, the completed tensor
algebra of $V$ is defined 
as 
$$ 
\oplus_{i\ge 0}V^{\hat\otimes n}, 
$$ 
and endowed with the obvious product, and denoted by $T(V)\hat{}$. 
These
objects are independent of the basis of neighborhoods chosen. 

$V\hat\otimes W$ (resp. $T(V)\hat{}$) is a separated complete
topological vector space (resp. algebra), with a basis of neighborhoods
of zero given by $\lim_{\leftarrow a,b}(V_{n}\otimes  W_{m}
/ V_{a} \otimes W_{b})$ (resp. the subalgebras generated by
the $V_{n}$). 

We will also define $V\bar\otimes W$ as the inverse limit 
$$
\lim_{\leftarrow a,b} 
(V \otimes W)/(V_{a} \otimes W + V\otimes W_{b}),
$$
whose topology is defined by 
the basis of neighborhoods
of zero given by $\lim_{\leftarrow a,b}(V_{n} \otimes W + V \otimes  W_{m})
/ (V_{a} \otimes W + V \otimes  W_{b})$. 

\subsection{Manin pairs and triples} \label{mp}

Let $X$ be a smooth, connected, compact complex curve, and $\omega$ be a 
nonzero
meromorphic differential on $X$. Let $S$ be a finite set of points of
$X$, containing the set $S_{0}$ of its zeros and
poles. For each $s\in S$, let $k_{s}$ be the local field at $s$ and
$$
k=\oplus_{s\in S}k_{s}. 
$$
Let $R$ be the ring of meromorphic functions on $X$, regular outside
$S$; $R$ can be viewed as a subring of $k$. $R$ is endowed with the
discrete topology and $k$ with its usual (formal series) topology. 
Let us define on $k$ the
bilinear form 
$$
\langle f,g \rangle_{k}=\sum_{s\in S}\res_{s}(fg\omega), 
$$
and the derivation 
$$
\pa f= df/\omega.
$$ We will use the notation ${\mathfrak x}(A)={\mathfrak x}\otimes A$,
for any ring $A$ over $\CC$ and complex Lie algebra $\mathfrak x$.

Let $\A={\mathfrak sl}_{2}(\CC)$. Define on $\A(k)$ the bilinear form
$\langle, \rangle_{\A(k)}$ by 
$$ 
\langle x\otimes \eps, y\otimes \eta\rangle_{\A(k)}=\langle
x,y\rangle_{\A}\langle \eps, \eta \rangle_{k}
$$
for $x,y\in\A, \eps,\eta\in k$, $\langle, \rangle_{\A}$ being the
Killing form of $\A$, the derivation $\pa_{\A(k)}$ 
by $\pa_{\A(k)}(x\otimes \eps)=x\otimes
\pa\eps$, for $x\in \A, \eps\in k$, 
and the cocycle
$$
c(\xi,\eta)=\langle \xi,\pa_{\A(k)} \eta\rangle_{\A(k)}. 
$$
Let $\hat{\G}$ be the central extension of $\A(k)$ by this cocycle. We
then have 
$$
\hat{\G}=\A(k)\oplus\CC K, 
$$
with bracket such that $K$ is central, and  $[\xi,\eta]
=([\bar \xi,\bar \eta],c(\bar \xi,\bar \eta)K)$, for any $\xi,\eta\in 
\hat{\G}$ with
first components $\bar \xi,\bar \eta$. 

Let us denote by $\pa_{\hat{\G}}$ the derivation of 
$\hat{\G}$ defined by $\pa_{\hat{\G}}(\xi,0)=(\pa_{\A(k)}\xi,0)$ and
$\pa_{\hat{\G}}(K)=0$. 

Let $\G$ be the skew product of $\hat{\G}$ with $\pa_{\hat{\G}}$. We
have 
$$
\G=\hat{\G}\oplus \CC D, 
$$
with bracket such that $\hat{\G}\to \G$, $\xi\mapsto (\xi,0)$ is a Lie
algebra morphism, and $[D,(\xi,0)]=(\pa_{\hat{\G}}(\xi),0)$ for
$\xi\in\hat{\G}$. 

View $\A(k)$ as a subspace of $\G=\hat{\G}\oplus \CC D=\A(k)\oplus \CC K
\oplus \CC D$, by $\xi\mapsto (\xi,0,0)$. 
Define on $\G$ the pairing $\langle, \rangle_{\G}$ by $\langle
K,D\rangle_{\G} =1$, $\langle K, \A(k)\rangle_{\G} = \langle D,
\A(k)\rangle_{\G} 
=0$, $\langle \xi,\eta \rangle_{\G}=\langle \xi,\eta \rangle_{\A(k)}$
for $\xi,\eta\in \A(k)$. 

Endow $\A(k)$ with $\langle, \rangle_{\A(k)}$. The subspace
$\A(R)\subset \A(k)$
is a maximal isotropic subalgebra of $\A(k)$, as follows from
Lemma \ref{statement-R}. 
Drinfeld's Manin pair is $(\A(k), \A(R))$ (see \cite{D}). 
In \cite{ER}, we introduced the following extension of this pair. 
Let $\G_{R}=\A(R)\oplus\CC D$; $\G_{R}\subset \G$ is a maximal isotropic
subalgebra of $\G$. The extended Drinfeld's Manin pair of 
\cite{ER} is then $(\G,\G_{R})$. 

In \cite{ER}, we also introduced the following
Manin triple. Let $\La$ be a Lagrangian complement to $R$ in $k$,
commensurable with $\oplus_{s\in S}\cO_{s}$ (where $\cO_{s}$ is the
completed local ring at $s$). 
Let $\N_{+}=\CC
e$, $\N_{-}=\CC f$, $\HH=\CC h$. Let 
$$
\G_{+}=\HH(R)\oplus \N_{+}(k)\oplus \CC D, 
\quad \G_{-}=(\HH\otimes \La)\oplus \N_{-}(k) \oplus \CC K,
$$
then $\G=\G_{+}\oplus \G_{-}$, and 
both $\G_{+}$ and $\G_{-}$ are maximal isotropic subalgebras of 
$\G$. The Manin triple is then $(\G,\G_{+},\G_{-})$. 

We will also consider the following Manin triple, that we may consider
as being obtained from the previous one by the action of the nontrivial
element of the Weyl group of $\A$. 
Let 
$$
\bar{\G}_{+}=\HH(R)\oplus \N_{-}(k)\oplus \CC D, 
\quad \bar{\G}_{-}=(\HH\otimes \La)\oplus \N_{+}(k) \oplus \CC K,
$$
then $(\G,\bar{\G}_{+},\bar{\G}_{-})$ again forms a Manin triple. 

\begin{remark} \label{frob}
{\it Generalizations.}
It is straightforward to generalize the centerless versions of 
the Manin pairs and triples
introduced above, as well as (as we will see) of their quantizations, 
to the case of a Frobenius algebra (i.e. a
commutative ring $k_{0}$ with a linear form $\theta\in (k_{0})^{*}$,
such that $(a,b)\mapsto \langle a,b \rangle_{k_{0}} = \theta(ab)$ is a
non-degenerate inner product), with a maximal isotropic subalgebra
$R_{0}$. 

It is also easy to generalize the Manin pairs and triples
$(\G,\G_{+},\G_{-}), (\G,\bar{\G}_{+}, \bar{\G}_{-})$ and $(\G,\G_{R})$,
as well as their quantizations in the sense of formal series, 
to the case where the Frobenius algebra is endowed with a derivation
$\pa_{0}$, such that $\theta\circ \pa_{0}=0$. \qed \medskip
\end{remark}

\subsection{Classical twists}

According to \cite{QG}, to each of the Manin triples
$(\G,\G_{+},\G_{-})$ and $(\G,\bar{\G}_{+},\bar{\G}_{-})$ is
associated a Lie bialgebra structure on $\G$; denote by
$\delta,\bar\delta:\G\to \G\hat\otimes\G$ the corresponding cocycle maps.

Let $\G_{\La} = (\A\otimes \La)\oplus \CC K \subset \G$; $\G_{\La}$ is a
Lagrangian complement of $\G_{R}$ in $\G$. It induces a Lie
quasi-bialgebra structure on $\G_{R}$, and 
from \cite{Kosmy} follows also that there is a Lie quasi-bialgebra
structure on $\G$, associated to the Manin pair $(\G,\G_{R})$ and to
$\G_{\La}$; we denote
by $\delta_{R}:\G\to\G\hat\otimes\G$ the corresponding cocycle map. 

These Lie (quasi-)bialgebra structures on $\G$ are related by the
following classical twist operations. 

Let $(e^{i})_{i\in \NN}, (e_{i})_{i\in \NN}$ be dual bases of $R$ and
$\La$; we choose them is such a way that $e_{i}$ tends to $0$ when $i$
tends to $\infty$.
Let $\eps^{i}, \eps_{i}, i\in \ZZ$ be dual bases of $k$, defined by
$\eps_{i}=e_{i}, \eps^{i}=e^{i}, i\ge 0$,  $\eps_{i}=e^{-i-1},
\eps^{i}=e_{-i-1}, i< 0$.

\begin{lemma} \label{class-twist}
Let $f=\sum_{i\in\ZZ}e[\eps^{i}]\otimes f[\eps_{i}]$; 
$f=f_{1}+f_{2}$, with 
$$
f_{1}=\sum_{i\in\NN}e[e_{i}]\otimes f[e^{i}], 
$$ 
and 
$$
f_{2}=\sum_{i\in\NN}e[e^{i}]\otimes f[e_{i}]. 
$$ 
For $\xi\in\G$, we have
$$
\delta_{R}(\xi)=\delta(\xi)+[f_{1}, \xi\otimes 1 + 1\otimes \xi], 
\quad 
\bar\delta(\xi)=\delta_{R}(\xi)+[f_{2}, \xi\otimes 1 + 1\otimes \xi]. 
$$
\end{lemma}

{\em Proof.} This is a consequence of the fact that 
the cocycle maps $\de,\bar\de,\de_{R}$ are respectively equal
to $\xi\mapsto [\sum_{j}\la_{j}\otimes \mu^{j},\xi \otimes 1 +
1\otimes \xi ]$, $ \xi \mapsto [\sum_{j}\la'_{j}\otimes \mu^{\prime j},\xi
\otimes 1 +
1\otimes \xi ]$, $ \xi\mapsto 
[\sum_{j}\la_{j}^{(R)}\otimes \mu^{j(R)},\xi \otimes 1 +
1\otimes \xi ]$, where $(\la_{j}),(\mu^{j})$,
resp. $(\la'_{j}),(\mu^{\prime j})$; $(\la_{j}^{(R)}), (\mu^{j(R)})$
are dual bases of $\G_{+}, \G_{-}$, resp. $\G'_{+}, \G'_{-}$; $\G_{R},
\G_{\La}$.   
\qed

\medskip\section{Quantization of Manin triples} \label{qmt}

\subsection{Results on kernels}

Recall that we have introduced dual bases  
$(e^{i})_{i\in \NN}, (e_{i})_{i\in \NN}$ of $R$ and
$\La$. 
Let $a_{0}\in R\hat\otimes \La$ be equal to 
$$
a_{0}=\sum_{i}e^{i}\otimes e_{i}.
$$
Note that $R\hat\otimes k$ is an algebra, to which belongs $a_{0}$.
Let 
$$
\gamma=(\pa\otimes 1)a_{0}-(a_{0})^{2}; 
$$
then 
\begin{lemma} (see \cite{ER}) \label{gamma-in-R}
$\gamma$ belongs to $R\otimes R$. 
\end{lemma}

Let $\hbar$ be a formal variable and let $T:k[[\hbar]]\to k[[\hbar]]$
be the operator equal to 
$$
T={{\sh(\hbar \partial)} \over {\hbar \partial}}.
$$
Let us use the notation $x\mapsto \wt{x}$ for the operation of
exchanging the two factors of $(k\bar\otimes k)[[\hbar]]$. 

\begin{prop}(see \cite{ER} and \cite{Some}, Prop. 1.11) \label{kernels}
For certain elements $\phi\in (R\otimes R)[[\hbar]]$, 
$\psi_{+},\psi_{-}\in \hbar (R\otimes R)[[\hbar]]$ depending on
$\g, \wt{\g}$ and their derivatives by universal formulas, we have the
following identities in $(R\hat\otimes k)[[\hbar]]$
$$ 
\sum_{i} Te^{i}\otimes e_{i}=
\phi+{1\over{2\hbar}}\ln{{1+a_{0}\psi_{-}} \over{1+a_{0}\psi_{+}}},
\quad \sum_{i} e^{i}\otimes Te_{i}=
-\wt{\phi}+{1\over{2\hbar}}\ln{{1-a_{0}\wt{\psi}_{+}}
  \over{1-a_{0}\wt{\psi}_{-}}} .
$$
\end{prop}

\begin{lemma} (see \cite{ER}) The expression 
$\sum_{i}Te^{i}\otimes e_{i} - e^{i}\otimes Te_{i}$
belongs to $S^{2}(R)[[\hbar]]$. We will denote by $\tau$ any element
of $(R\otimes R)[[\hbar]]$, such that 
\begin{equation} \label{id-tau}
\tau+\wt{\tau}=\sum_{i}Te^{i}\otimes e_{i} - e^{i}\otimes Te_{i},
\end{equation}
and define the linear map $U:\La[[\hbar]] \to R[[\hbar]]$ by 
\begin{equation} \label{U}
U\la=\langle\tau, 1\otimes \la \rangle. 
\end{equation}
\end{lemma}

Note that $\sum_{i}Te^{i}\otimes e_{i}$ is well-defined in
$(R\hat\otimes k)[[\hbar]]$, because $e_{i}$ tends to zero as $i$ tends
to infinity. Since $\pa$ is a continuous map from $k$ to itself, the
same is true for the sequence $\pa^{k}e_{i}$. So $\sum_{i}e^{i}\otimes
Te_{i}$ is well-defined in the same space; 
$\sum_{i\in \ZZ}T\eps^{i}\otimes \eps_{i} -\eps^{i}\otimes
T\eps_{i}$
is well-defined in $(k\bar\otimes k)[[\hbar]]$ for the same reasons. 

\begin{remark} (see \cite{ER})
Let $q=e^{\hbar}$, and 
$f$ be the function defined by $f(x)={{q^{\pa}-1}\over
{\hbar \pa}}$ and 
$\tau_{0}=\sum_{i}f(\pa)e^{i}\otimes e_{i}-e^{i}\otimes f(-\pa)e_{i}$ (since
$\pa$ is a continuous map from $k$ to itself, each $\sum_{i}e^{i}\otimes
\pa^{k}e_{i}$ is well-defined in $R\hat\otimes k$, so that $\tau_{0}$ is
well-defined in $(R\hat\otimes k)[[\hbar]]$). 

Note that the
formal series $f(\pa_{z})-f(-\pa_{w})$ is divisible by $\pa_{z}+\pa_{w}$
in $\CC[\pa_{z},\pa_{w}][[\hbar]]$, and denote their ratio by  
${{f(\pa_{z})-f(-\pa_{w})}\over{\pa_{z}+\pa_{w}}}$. 

Attach indices $z$ and $w$ to the first and second factors of $(k\hat\otimes
k)[[\hbar]]$. We have
$$
\tau_{0}={{f(\pa_{z})-f(-\pa_{w})}\over{\pa_{z}+\pa_{w}}}
(\gamma-\wt{\gamma}), 
$$
and $\tau_{0}$ satisfies the identity (\ref{id-tau}). 
\qed\medskip
\end{remark}

Consider now the quantity $\exp(2\hbar \sum_{i\in\NN}
e^{i}\otimes (T+U)e_{i})$; it
belongs to $(R\hat\otimes k)[[\hbar]]$. 
Let for each $s\in S$, $z_{s}$ be a local coordinate on $X$ near $s$. 

\begin{lemma} \label{regularisation}
For any $\al\in k$, we have
\begin{align} \label{reg}
  (\al\otimes 1 -1 \otimes \al +\psi_{-}(\al\otimes 1 & -1 \otimes
  \al) a_{0}) q^{2\sum_{i}(T+U)e_{i}\otimes e^{i}} \\ \nonumber &
  =(\al\otimes 1- 1 \otimes \al +\psi_{+}(\al \otimes 1 -1 \otimes
  \al)a_{0})q^{2(\tau-\phi)}.
\end{align}

For any $\al\in k$, $(\al \otimes 1 -1 \otimes \al)a_{0}$ 
belongs to $\prod_{s,t\in S}\CC((z_{s},w_{t}))$. If $\al$ belongs to $R$, then
$(\al \otimes 1 -1 \otimes \al)a_{0}$ belongs to $R\otimes R$. 
\end{lemma}

{\em Proof.} The first part of the lemma is proved using
that $\tau=\sum_{i\in\NN}Ue_{i}\otimes e^{i}$, the second identity of
Prop. \ref{kernels} and 
\begin{equation}
\label{delta}
(\al\otimes 1 - 1 \otimes \al)\wt{a}_{0}
=
-(\al\otimes 1 - 1 \otimes \al)a_{0}, \quad \forall \al\in k
\end{equation}
((\ref{delta}) follows from the fact that $a_{0}+\wt a_{0}$ verifies
$\langle a_{0}+\wt a_{0}, id \otimes \beta \rangle_{k}=\beta$ for any
$\beta\in k$; the product is taken with the second components of a
decomposition of $a_{0}+ \wt a_{0}$). 

Let us pass to the second part of the lemma. For $N$ integer, set
$$
k_{N}=\prod_{s\in S}z_{s}^{-N}\CC[[z_{s}]].
$$ 
For any $\al\in k$, there
exists an integer $N_{\al}$ such that 
$(\al\otimes 1 - 1 \otimes \al)a_{0}$ belongs to $k\hat\otimes
k_{N_{\al}}$. Since $(\al\otimes 1 - 1 \otimes \al)\wt{a}_{0}$ belongs
to $k_{N_{\al}}\bar\otimes k$, and using (\ref{delta}), we obtain 
that 
$(\al\otimes 1 - 1 \otimes \al)\wt{a}_{0}$ belongs to 
$ ( k\bar\otimes
k_{N_{\al}} ) \cap ( k_{N_{\al}}\bar\otimes k ) 
\subset \prod_{s,t\in
S}\CC((z_{s},w_{t}))$. 

The third part of the lemma follows from the following statement. 
Let $\al\in R$, then for any $\beta\in R$, one checks that 
$$
\langle (\al\otimes 1 - 1 \otimes \al)a_{0}, \beta \otimes id
\rangle_{k}
=
\langle (\al\otimes 1 - 1 \otimes \al)a_{0}, id \otimes \beta
\rangle_{k}=0,
$$
where the products are taken with the first and second components of a
decomposition of $(\al\otimes 1 - 1\otimes \al)a_{0}$. 
\qed\medskip

\begin{remark} \label{terminology}
Attach indices $z$ and $w$ to the first and second factors of
$k\bar\otimes k$. Let $P$ be a differential operator, acting on $k$.
The quantity $K_{P}(z,w)=\sum_{i\in\ZZ}P\eps^{i}\otimes
\eps_{i}=\sum_{i\in\ZZ}(P\eps^{i})(z)\eps_{i}(w)$ can be considered as a
kernel for the operator $P$, because of the identity
$$
(Pf)(w)=\res_{w\in S} K_{P}(z,w)f(w)\omega_{w}, \quad \forall f\in k.
$$
Suppose that $P$ preserves $R$. Then 
 $\bar K_{P}(z,w)=\sum_{i\in\ZZ}Pe^{i}\otimes
e_{i}=\sum_{i\in\ZZ}(Pe^{i})(z)e_{i}(w)$ can also be considered as a
kernel for $P$, because 
$$
(Pf)(w)=\res_{w\in S} \bar K_{P}(z,w)f(w)\omega_{w}, \quad \forall f\in R.
$$
\qed\medskip 
\end{remark}

\begin{remark} \label{Hochschild}
The map $R\mapsto R \otimes R$, $\al\mapsto (\al\otimes 1 - 1
\otimes \al)a_{0}$ defines a nontrivial element of the Hochschild cohomology
$H^{1}(R,R\otimes R)$ (where $R\otimes R$ is the $R$-bimodule, with
left action by multiplication on the first factor on the tensor
product 
and right action by multiplication on the second one). 

Note that there is a natural map $H^{1}(R, R\otimes R)\to H^{1}(k,
k\hat\otimes_{R}(R\otimes R)\hat\otimes_{R}k)$; $k\hat\otimes_{R}(R\otimes
R)\hat\otimes_{R}k$ is equal to $\prod_{s,t\in S}\CC((z_{s},w_{t}))$, so the
image by this map of the above class is nonzero
[since $a_{0}\notin \prod_{s,t\in S}\CC((z_{s},w_{t}))$]. 
\end{remark}

\subsection{Presentation of $U_{\hbar}\G$}

In \cite{ER}, we introduced a Hopf algebra $U_{\hbar}\G$ quantizing
$(\G,\delta)$. 

It is the quotient of $T(\G)\hat{}[[\hbar]]$ by the following
relations. Let $e,f,h$ be the Chevalley basis of ${\mathfrak sl}_{2}(\CC)$.  
Denote in $T(\G)\hat{}[[\hbar]]$, the element $x\otimes \eps\in \A(k)\subset
\G$ of
$\G$ by $x[\eps]$ and let for $r\in R$, $h^{+}[r]=h[r]$,
$h^{-}[\la]=h[\la]$. 
Introduce the generating series 
$$
e(z)=\sum_{i\in\ZZ}e[\eps_{i}]\eps^{i}(z), \quad
f(z)=\sum_{i\in \ZZ}f[\eps_{i}]\eps^{i}(z), 
$$
$$
h^{+}(z)=\sum_{i\in\NN}h^{+}[e^{i}]e_{i}(z), \quad
h^{-}(z)=\sum_{i\in\NN}h^{-}[e_{i}]e^{i}(z).
$$
The relations for $U_{\hbar}\G$ are (Fourier modes of)
\begin{equation} \label{comm}
[h^{+}[r], h^{+}[r']]=0, 
\end{equation}
\begin{equation} \label{h-h}
[K,\on{anything}]=0, \quad
[h^{+}[r], h^{-}[\la]]= {2 \over  \hbar}   \langle 
(1 - q^{- K  \partial })r,  \la  \rangle, 
 \end{equation}
\begin{equation} \label{h-h-}
[h^{-}[\la], h^{-}[\la']]={2\over\hbar} \left(
\langle
T((q^{K\pa}\la)_{R}), q^{K\pa}\la' 
\rangle
+\langle U \la, \la' \rangle
-\langle U ( (q^{K\pa}\la)_{\La} ) , q^{K\pa}\la' \rangle
\right), 
\end{equation}
\begin{equation}   \label{h-e}
[h^{+}[r], e(w)]=2 r(w)e(w), \quad [h^{-}[\la],
e(w)]=2[(T+U)(q^{K\pa}\la)_{\La}](w)e(w), 
\end{equation}
\begin{equation}   \label{h-f}
[h^{+}[r], f(w)]=-2 r(w)f(w), \quad [h^{-}[\la],
f(w)]=-2[(T+U)\la](w)f(w), 
\end{equation}
\begin{align}  \label{e-e}
[z_{s}-w_{s} & +\psi_{-}(z,w)a_{0}(z,w)(z_{s}-w_{s})]
 e(z)e(w)
\\ \nonumber 
& =q^{2(\tau-\phi)(z,w)}
[z_{s}-w_{s}+\psi_{+}(z,w)a_{0}(z,w)(z_{s}-w_{s})]
e(w)e(z), \quad \forall s\in S 
\end{align}
\begin{align}  \label{f-f}
q^{2(\tau-\phi)(z,w)}
[z_{s}-w_{s} & +\psi_{+}(z,w)a_{0}(z,w)(z_{s}-w_{s})]
f(z)f(w)
\\ \nonumber 
&=
[z_{s}-w_{s} +\psi_{-}(z,w)a_{0}(z,w)(z_{s}-w_{s})]
 f(w)f(z), \quad \forall s\in S 
\end{align}
\begin{equation}   \label{e-f}
[e(z),f(w)]={1\over \hbar} [ \delta(z,w)q^{((T+U)h^{+})(z)}
-(q^{-K\pa_{w}}\delta(z,w))q^{-h^{-}(w)} ] ,
\end{equation}
\begin{equation} \label{D-h}
[D,h^{+}[r]]=h^{+}[\pa r], [D, h^{-}(z)]=-(\pa h^{-})(z)-
\sum_{i\in \NN}[(1+q^{-K\pa})Ae_{i}](z)h^{+}[e^{i}]+\cA(z), 
\end{equation}
\begin{equation} \label{D-ef}  
[D, x^{\pm}(z)]=-(\pa x^{\pm})(z)+{\hbar}\sum_{i\in
\NN}(Ae_{i})(z)h^{+}[e^{i}] x^{\pm}(z), 
\end{equation}
for any $r,r'\in R$, $\la,\la'\in \La$, $x^{\pm}=e,f$, 
with $A: \La[[\hbar]] \to
R[[\hbar]]$ defined by 
\begin{equation} \label{A}
A\la = T((\pa \la)_{R})+\pa (U \la) - U((\pa \la)_{\La}), \quad \forall
\la\in \La, 
\end{equation}
and 
\begin{equation} \label{cal-A}
\cA(z)= \sum_{i\in \NN}
[(1+q^{-K\pa})Ae_{i}](z)[(1-q^{-K\pa})e^{i}](z)
+{2\over \hbar}[(q^{-K\pa}-1)Ae_{i}](z)e^{i}(z). 
\end{equation}
Note that $A$ is anti-self-adjoint, so that $\sum_{i\in \NN}A
e^{i}\otimes e_{i} = - \sum_{i\in \NN} e^{i}\otimes Ae_{i}$, and
$\sum_{i\in \NN}(A e_{i})(z) e^{i}(z)=0$. 

For any positive integer $n$, 
define the completed tensor power $(U_{\hbar}\G)^{\hat\otimes n}$ as the
quotient of $T(\G^{n})\hat{}[[\hbar]]$ by the ideal generated by the usual $n$
copies of the above relations. 

The formulas 

\begin{equation}  \label{Delta-K}
\Delta(K)=K\otimes 1 + 1 \otimes K
\end{equation}
\begin{equation}    \label{Delta-h}
\Delta(h^{+}[r])=h^{+}[r]\otimes 1+1\otimes h^{+}[r], \quad
\Delta(h^{-}(z))=h^{-}(z)\otimes 1+1\otimes (q^{-K_{1}\pa}h^{-})(z), 
\end{equation}
\begin{equation}   \label{Delta-e}
\Delta(e(z))=e(z)\otimes q^{((T+U)h^{+})(z)} + 1\otimes e(z), 
\end{equation}
\begin{equation}   \label{Delta-f}
\Delta(f(z))=f(z)\otimes 1+ q^{-h^{-}(z)}\otimes (q^{-K_{1}\pa}f)(z),   
\end{equation}
\begin{equation} \label{Delta-D}
\Delta(D)= D \otimes 1 + 1 \otimes D + \sum_{i\in \NN}{\hbar \over 4}
h^{+}[e^{i}] \otimes h^{+}[A e_{i}],
\end{equation} 
$r\in R$, for the coproduct, 
\begin{equation} \label{counit}
\varepsilon(h^{+}[r])=\varepsilon(h^{-}[\la])=\varepsilon(x[\eps])=
\varepsilon(D)=\varepsilon(K)=0,
\end{equation}
$x=e,f$, $r\in R, \la\in\La, \eps\in k$, 
for the counit,
\begin{equation} \label{S-efhKD}
S(h^{+}[r])=-h^{+}[r],\ (Sh^{-})(z)=-(q^{K\pa}h^{-})(z), \
S(D)=-D + {\hbar \over 4}\sum_{i\in
\NN}h^{+}[e^{i}]h^{+}[Ae_{i}],
\end{equation} 
\begin{equation} \label{S-ef}
(Se)(z)=-e(z)q^{((T+U)h^{+}) (z)}, \quad
(Sf)(z)=-\left(q^{K\pa}(q^{h^{-}} f) \right)(z), 
\quad S(K)=-K,  
\end{equation} 
$r\in R$, 
for the antipode, 
define a topological (with respect to the completion introduced above)
Hopf algebra structure on $U_{\hbar}\G$.

\begin{notation} \label{notation}
We have posed $\delta(z,w)=\sum_{i\in \ZZ}\eps^{i}(z)\eps_{i}(w)$; this
is an element of $k\bar\otimes k$. 
The indices $R$ and $\La$ denote the projections on the first and second
factor of the decomposition $k[[\hbar]]=R[[\hbar]]\oplus \La[[\hbar]]$. 

In (\ref{e-e}), (\ref{f-f}), we have attached indices $z$ and $w$ to the
first and second factors of $k\bar\otimes k$. Recall that
$a_{0}(z,w)(z_{s}-w_{s})$ belongs to $\CC((z_{s},w_{s}))$. 

$K_{1},K_{2}$ respectively mean $K\otimes 1, 1\otimes K$. 
The $K$, $f(z)$ and $h^-(z)$ used here correspond to 
$2K$, ${1 \over \hbar} (q^{-K \partial}f)(z)$ and 
$(q^{-K \partial} h^-)(z)$ of \cite{ER} respectively. 
\qed\medskip
\end{notation}

\begin{remark} \label{other-e-e}
{\it Variants of the vertex relations (\ref{e-e}) and (\ref{f-f}). }
Due to the Hochschild cocycle properties explained in
rem. \ref{Hochschild}, 
relations (\ref{e-e}) and (\ref{f-f}) are equivalent to the following
ones, 
\begin{align}  \label{e-e-variant}
[ & \al(z)-\al(w) +\psi_{-}(z,w)a_{0}(z,w)(\al(z)-\al(w))]
 e(z)e(w)
\\ \nonumber 
& =q^{2(\tau-\phi)(z,w)}
[\al(z)-\al(w)+\psi_{+}(z,w)a_{0}(z,w)(\al(z)-\al(w))]
e(w)e(z), \quad \forall \al\in k
\end{align}
\begin{align}  \label{f-f-variant}
q^{2(\tau-\phi)(z,w)}
[ & \al(z)-\al(w) +\psi_{+}(z,w)a_{0}(z,w)(\al(z)-\al(w))]
f(z)f(w)
\\ \nonumber 
&=
[\al(z)-\al(w) +\psi_{-}(z,w)a_{0}(z,w)(\al(z)-\al(w))]
 f(w)f(z), \quad \forall \al\in k; 
\end{align}
note that for any $\al\in k$, $a_{0}(z,w)(\al(z)-\al(w))$  
belongs to $\prod_{s,t\in S}\CC((z_{s},w_{t})) = k \hat\otimes k$. 
\qed\medskip
\end{remark}

\begin{remark} \label{informal-writing-e-e}
Due to (\ref{reg}), the $e-e$ and $f-f$ relations
(\ref{e-e}) and (\ref{f-f}) can be informally written as 
\begin{align}  \label{informal-vx}
e(z)e(w) =q^{2\sum_{i}((T+U)e_{i})(z)e^{i}(w)}e(w)e(z), 
f(z)f(w)=q^{2\sum_{i}e^{i}(z)((T+U)e_{i})(w)}f(w)f(z). 
\end{align}
\qed\medskip
\end{remark}

\begin{remark} \label{rewriting-D-ef}
The relations (\ref{D-ef}) can be rewritten as 
$$
[D,q^{-((T+U)h^{+})(z)} x(z)] = - \pa_{z} ( q^{-((T+U)h^{+})(z)} x(z) ), 
$$
$x=e,f$. 
\end{remark}

\subsection{Presentation of $U_{\hbar}\bar{\G}$} \label{ug}

The Lie bialgebra $(\G, \bar\delta)$ also admits a quantization. We denote
by $U_{\hbar}\bar{\G}$ the corresponding Hopf algebra. It is the quotient of
$T(\G)\hat{}[[\hbar]]$ by the following relations. Let us overline in
the case of $U_{\hbar}\bar{\G}$, the
notation used in the case of $U_{\hbar}\G$. 
The algebra relations are
identical to those of $U_{\hbar}\G$, except for those involving $\bar
h^{-}[\la]$; we have 
\begin{equation} \label{bar-h-e} 
[\bar h^{-}[\la], \bar e(z)] = 2 ((T+U)\la)(z)\bar e(z),
\end{equation}
\begin{equation} \label{bar-h-f} 
[\bar h^{-}[\la], \bar f(z)] =  - 2 ((T+U)(q^{-\bar K
\pa}\la)_{\La})(z)\bar f(z),
\end{equation}
\begin{equation} \label{bar-h-h+} 
[\bar h^{+}[r], \bar h^{-}[\la]] =  { 2 \over \hbar }
\langle  (q^{\bar K\pa}-1)r , \la \rangle, 
\end{equation}
\begin{equation} \label{bar-h-h-}
[\bar h^{-}[\la], \bar h^{-}[\la']]={2\over\hbar} \left(
\langle U ( (q^{- \bar K\pa}\la)_{\La} ) , q^{ - \bar K\pa}\la' \rangle
- \langle U \la, \la' \rangle
- \langle
T((q^{ - \bar K\pa}\la)_{R}), q^{ - \bar K\pa}\la' 
\rangle
\right), 
\end{equation}
\begin{equation} \label{bar-e-f}
[\bar e (z), \bar f (w) ]={1\over \hbar} [\delta(z,w)
q^{((T+U)\bar h^{+})(z) } - 
q^{-\bar K \pa_{w}}\left( q^{-\bar h^{-}(w)}\de(z,w)\right),
\end{equation}
\begin{equation} \label{bar-D-h}
[\bar D,\bar h^{-}(z)] = - (\pa\bar h^{-}) (z) - [(1+q^{\bar K
\pa}) A e_{i}](z) \bar h^{+} [e^{i}] + \bar{\cal A}(z), 
\end{equation}
where
\begin{equation} \label{barcA}
\bar{\cal A}(z)= \sum_{i\in\NN}
[(1+q^{\bar K \pa})Ae_{i}](z)[(q^{\bar K
\pa}-1)e_{i}](z) + {2\over \hbar}[(1-q^{\bar K\pa})Ae_{i}](z) (q^{\bar
K\pa}e^{i})(z).
\end{equation}
The coalgebra structure of $U_{\hbar}\bar{\G}$  is defined by the
coproduct 
\begin{equation} \label{Delta'-K}
\bar\Delta(\bar K)=\bar K\otimes 1 + 1 \otimes \bar K, 
\end{equation}
\begin{equation}    \label{Delta'-h}
\bar\Delta(\bar h^{+}[r])=\bar h^{+}[r]\otimes 1+1\otimes \bar h^{+}[r], \quad
\bar\Delta(\bar h^{-}(z))=(q^{\bar K_{2}\pa}\bar h^{-})(z)\otimes
1+1\otimes \bar h^{-}(z), 
\end{equation}
\begin{equation}   \label{Delta'-e}
\bar\Delta(\bar e(z))=(q^{\bar K_{2}\pa}\bar e)(z)\otimes q^{-\bar
h^{-}(z)} + 1\otimes \bar e(z), 
\end{equation}
\begin{equation}   \label{Delta'-f}
\bar\Delta(\bar f(z))=\bar f(z)\otimes 1+ q^{((T+U)\bar h^{+})(z)} \otimes
\bar f(z),
\end{equation}
\begin{equation} \label{Delta'-D}
\bar\Delta(\bar D)= \bar D \otimes 1 + 1 \otimes \bar D 
+ \sum_{i\in \NN}{\hbar \over 4}
\bar h^{+}[e^{i}] \otimes \bar h^{+}[A e_{i}],
\end{equation} 
$r\in R$, the counit
\begin{equation} \label{counit'}
\bar\varepsilon(\bar h^{+}[r])=\bar \varepsilon(\bar h^{-}[\la])=\bar
\varepsilon(\bar x[\eps])=
\bar \varepsilon(\bar D)=\bar \varepsilon(\bar K)=0,
\end{equation}
$x=e,f$, $r\in R, \la\in\La, \eps\in k$,
and the antipode 
\begin{equation} \label{S'-hKD}
\bar S(\bar h^{+}[r])=-\bar h^{+}[r],\ (\bar S\bar h^{-})(z)=-(q^{-\bar
K\pa}\bar h^{-})(z), \
\bar S(\bar D)=-\bar D + {\hbar \over 4}\sum_{i\in
\NN}\bar h^{+}[e^{i}]\bar h^{+}[Ae_{i}],
\end{equation} 
\begin{equation} \label{S'-ef}
(\bar S\bar e)(z)= \left( q^{-\bar K \pa} (\bar e q^{\bar h^{-}}) \right)
(z), 
\quad
(\bar S\bar f)(z) = - q^{-((T+U)\bar h^{+})(z)}\bar f (z), 
\quad \bar S(\bar K)=-\bar K,  
\end{equation} 
$r\in R$. 

Then 

\begin{prop} (see \cite{ER})
  Define $U_{\hbar}\G_{\pm}$ and $U_{\hbar}\bar\G_{\pm}$ as the
  subalgebras of $U_{\hbar}\G$ and $U_{\hbar}\bar\G$ generated by
  $\G_\pm$ and $\bar\G_\pm$, respectively. The pairs
  $(U_{\hbar}\G_+,\Delta)$ and $(U_{\hbar}\G_-,\Delta')$, as well as
  $(U_{\hbar}\bar\G_+,\bar\Delta)$ and
  $(U_{\hbar}\bar\G_-,\bar\Delta')$, form dual Hopf algebras,
  quantizing the Lie bialgebra structures defined by
  $(\G_{\pm},\pm\delta)$ and $(\bar\G_{\pm},\pm\bar\delta)$. The pairs
  $(U_\hbar\G,\Delta)$ and $(U_\hbar\bar\G,\Delta)$ are the double
  Hopf algebras of $(U_{\hbar}\G_+,\Delta)$ and
  $(U_{\hbar}\bar\G_+,\bar\Delta)$ respectively, and define
  quantizations of the Lie bialgebras $(\G,\delta)$ and
  $(\bar\G,\bar\delta)$.
\end{prop}

Here $\Delta'$ denotes $\Delta$ composed with the permutation of
factors.

Then 

\begin{prop} \label{isom}
The map 
$$
x[\eps]\mapsto \bar x[\eps], h^{+}[r]\mapsto \bar h^{+}[r],
K \mapsto \bar K, D\mapsto \bar D, h^{-}[\la]\mapsto \bar
h^{-}[(q^{K\pa}\la)_{\La}],
$$ 
$x=e,f$, $\eps\in k, r\in R,\la\in\La$,
extends to an algebra isomorphism from $U_{\hbar}\G$ to $U_{\hbar}\bar{\G}$.
\end{prop}

In what follows, we will denote elements of $U_{\hbar}\bar{\G}$
as elements of $U_{\hbar}\G$, implicitly making use of this
isomorphism. 


\medskip\section{PBW results for $U_{\hbar}\G$}

\subsection{PBW result for algebras presented by vertex relations}
\label{vx} 

We will now prove a PBW statement, which was used implicitly in
\cite{ER}. 

Let $\zeta$ be an indeterminate, and let $V$ be the field of Laurent
series $\CC((\zeta))$. Let $\gamma_{n}=\zeta^{n}, n\in \ZZ$,
and
let us organize the $(\g_{n})_{n\in \ZZ}$ in the generating series 
$$
\g(z)=\sum_{i\in \ZZ}\g_{n}z^{-n}.
$$
Let $\hbar$ be a formal variable, and $\cA$ be the quotient of
$T(V)\hat{}[[\hbar]]$ by the relations obtained as the Fourier modes of
\begin{equation}  \label{relations}
(z-w+A(z,w))\g(z)\g(w)=(z-w+B(z,w))\g(w)\g(z), 
\end{equation}
for $A,B\in \hbar \CC((z,w))[[\hbar]]$. 

\begin{lemma} \label{generating}
Assume that $A,B$ satisfy the relation 
\begin{equation} \label{compat}
(z-w+B(z,w))(w-z+B(w,z))=(z-w+A(z,w))(w-z+ A(w,z)), 
\end{equation}
and the series $z-w+A$ and $z-w+B$ do not define the same ideal 
of $\CC((z,w))[[\hbar]]$. 
Then any element of $\cA$ can be written as a sum 
$$ \sum_{p=0}^{k} \sum_{i_{1}< \ldots< i_{p}, \al_{i}\ge 1}\la_{i_{1},
  \ldots, i_{p}}^{(\al_{1}, \ldots, \al_{p})}
\g_{i_{1}}^{\al_{1}}\ldots \g_{i_{p}}^{\al_{p}},
$$ $k\ge 0$, $\la_{i_{1}, \ldots, i_{p}}^{(\al_{1}, \ldots, \al_{p})}$
scalars, such that the number of indices $((i_{1},\ldots,
i_{p}),(\al_{1},\ldots, \al_{p}))$ with $i_{1}=M$ and $\la_{i_{1},
  \ldots, i_{p}}^{(\al_{1}, \ldots, \al_{p})}\neq 0$, is finite for
all $M$ and zero for $M$ large enough.
\end{lemma}

\noindent{\em Proof.} (\ref{compat}) implies that the divisors of 
$z-w+A$ and $z-w-\wt B$ coincide. Therefore, the relations
(\ref{relations}) can be put is the form
\begin{equation} \label{rel'}
\bar q_+(z,w) \gamma(z) \gamma(w) = \bar q_-(z,w) \gamma(w) \gamma(z), 
\end{equation}
with $\bar q_+(z,w) = z-w + \bar C(z,w)$ and $\bar q_- =  -\kappa \bar 
q_+^{(21)}$, 
with $\bar C\in \hbar\CC((z,w))[[\hbar]]$ and $\kappa\in 1 + \hbar
\CC((z,w))[[\hbar]]$, $\kappa \kappa^{(21)} = 1$.

Set $q_+ = \kappa^{-1/2} \bar q_+$, $q_- = \kappa^{-1/2} \bar q_-$, then 
relations (\ref{rel'}) can be written as 
\begin{equation} \label{rel''}
q_+(z,w) \gamma(z) \gamma(w) = q_-(z,w) \gamma(w) \gamma(z), 
\end{equation}
with $q_-(z,w) = - q_+(w,z)$. 

Let us prove now that the $\g_{i_{1}}^{\al_{1}}\ldots \g_{i_{p}}^{\al_{p}}$
form a generating family of $\cA$. 

Let $\cA_{2}$ be the span in $\cA$ of infinite
series $\sum_{n,m\ge N}a_{pq}e_{p}e_{q}$. 
(\ref{relations}) allows us to write, for any $n,m\in \ZZ$, 
$$
\rho_{n,m}=[\g_{n+1},\g_{m}]-[\g_{n},\g_{m+1}]
$$
as a series in $\hbar \cA_{2}[[\hbar]]$. 
We rearrange this system of relations in the following way. Let
$$
\tau_{n-k,n+k}=\rho_{n-k,n+k}+\rho_{n-k+1,n+k-1}+\ldots+\rho_{n-1,n+1} 
$$ for $k>0$,
$$
\tau_{n+1-k,n+k}=\rho_{n+1-k,n+k}+\rho_{n+2-k,n+k-1}+\ldots+\rho_{n,n+1}
$$
for $k\ge 1$, 
$$ \tau_{n+k,n-k}=\rho_{n+k,n-k}+\rho_{n+k-1,n-k+1}+\ldots+\rho_{n,n}
$$
for $k\ge 0$,  
$$
\tau_{n+1+k,n-k}=\rho_{n+1+k,n-k}+\rho_{n+k,n-k+1}+\ldots+\rho_{n+1,n}
$$
for $k\ge 0$, 
then the system of expressions for the $\rho_{n,m}$ is equivalent to as
a system of
expressions for the $\tau_{n,m}$.

Note that 
$$
\tau_{n-k,n+k}=-[\g_{n-k},\g_{n+k+1}]+[\g_{n},\g_{n+1}]
$$
for $k> 0$,  
$$
\tau_{n+1-k,n+k}=[\g_{n+1-k},\g_{n+k+1}]
$$
for $k\ge 1$, 
$$
\tau_{n+k,n-k}=[\g_{n+k+1},\g_{n-k}]-[\g_{n},\g_{n+1}]
$$
for $k\ge 0$,  
$$
\tau_{n+1+k,n-k}=[\g_{n+2+k}, \g_{n-k}]
$$
for $k\ge 0$. 

The expression for $\tau_{n,n}$ yields an expression for
$[\g_{n},\g_{n+1}]$. Substracting this expression to the expressions
for $\tau_{n-k,n+k}$ and adding to the expression for
$\tau_{n+k,n-k}$, we derive expressions for the $[\g_{n},\g_{m}]$.
This means that an arbitrary monomial in the $\g_{i}$'s can be
expressed as a linear combination of the $(\g_{i_{1}}^{\al_{1}}\ldots
\g_{i_{p}}^{\al_{p}})_{i_{1}<\ldots< i_{p},\al_{i}\ge 1}$ of the form
described. \hfill  \qed \medskip 

Let us now prove that $(\gamma_{i_1}^{\al_1} \cdots
\gamma_{i_p}^{\al_p})_{i_1 < i_2 < \cdots, \al_i \geq 1}$ forms a
topological basis of $\cA$. For this, we will construct an isomorphism of
$\cA$ with a Feigin-Odesski-type (of shuffle) algebra (see \cite{FO}).
Define $FO$ as the direct sum $\oplus_{n\ge 0} FO^{(n)}$, where
$FO^{(n)}$ is the the subspace of $\CC((z_1)) \cdots ((z_n))[[\hbar]]$
formed by the elements $\la$ such that $\prod_{i<j}q_-(z_i,z_j) \la$
belongs to $k^{\hat\otimes n}[[\hbar]]$ and is totally antisymmetric in
$z_1,\ldots,z_n$.

Define a product of $FO$ as follows: let $\eps$ and $\eta$ belong to
$FO^{(n)}$ and $FO^{(m)}$, then their product $\eps * \eta$ lies in
$FO^{(n+m)}$ and is equal to
\begin{align} \label{shuffle}
& (\eps * \eta)(z_1,\cdots,z_{n+m}) \\ & \nonumber ={{n!m!}\over{(n+m)!}} 
\sum_{\sigma\in Sh_{n,m}} \prod_{i<j,\si(i)>\si(j)} 
q(z_i,z_j)\eps(z_{\si(1)} , \cdots , z_{\si(n)}) 
\eta(z_{\si(n+1)} , \cdots , z_{\si(n+m)}) , 
\end{align}
where $Sh_{n,m}$ is the set of shuffle transformations of
$\{1,\ldots,n+m\}$, that is the set of permutations $\si$ of that set
such that $\si(i) < \si(j)$ is $i<j\leq n$ if $m+1\leq i<j$, and
$q(z,w)$ is the ratio $q_+(z,w)/q_-(z,w)$, expanded for $w<<z$.  It is
then clear that the product (\ref{shuffle}) is well-defined and
associative.

The space $FO^{(n)}$ can also be described as follows: 
\begin{prop} \label{str:FO}
  Define a symmetrization map $S_q$ from $k^{\hat\otimes n}[[\hbar]]$
  to $\CC((z_1)) \cdots ((z_n))[[\hbar]]$ by 
$$ 
S_{q}(\eps) = \sum_{\sigma\in
    \SS_n}\eps(z_{\sigma(1)},\cdots,z_{\sigma(n)}) \prod_{i<j, \sigma(i)
    > \sigma(j)}q(z_i,z_j)
$$ 
(recall that
  $k^{\hat\otimes n} =
  \CC[[z_1,\cdots,z_n]][z_1^{-1},\cdots,z_n^{-1}]$). 

  Then $FO^{(n)}$ is both equal to the image by $S_q$ of
  $k^{\hat\otimes n}[[\hbar]]$ and of its subspace of totally
  symmetric elements.
\end{prop}

{\em Proof.} It is clear that both images lie in $FO^{(n)}$. To prove
that any element of $FO^{(n)}$ is in $S_q(k^{\hat\otimes
  n}[[\hbar]]^{\SS_n})$, pick any totally antisymmetric $a$ in
$k^{\hat\otimes n}$; we should write it in the form
\begin{equation} \label{eqn} 
a = \left( \sum_{\si\in \SS_n} \prod_{i<j,\si(i)>\si(j)} q_{+}(z_i,z_j)
\prod_{i<j,\si(i)<\si(j)} q_{-}(z_i,z_j) \right) 
\eps(z_1,\cdots,z_n) , 
\end{equation}
with $\eps$ totally symmetric in $k^{\hat\otimes n}[[\hbar]]$. The sum
$$
\sum_{\si\in \SS_n} \prod_{i<j,\si(i)>\si(j)} q_{+}(z_i,z_j)
\prod_{i<j,\si(i)<\si(j)} q_{-}(z_i,z_j)
$$ lies in $n! \prod_{i<j}(z_i - z_j) + \hbar k^{\hat\otimes
  n}[[\hbar]]$. Since it is also totally antisymmetric, it is of
the form $n! \prod_{i<j}(z_i - z_j) (1 + \sum_{i\ge 1}\hbar^i
\kappa_i)$, $\kappa_i \in k^{\hat\otimes n}$. Since $a$ is totally
antisymmetric, we can write it as a product $\prod_{i<j} (z_i - z_j)
\cdot s$, $s$ totally symmetric in $k^{\hat\otimes n}[[\hbar]]$. We
then set $\eps = s / (1 + \sum_{i\ge 1}\hbar^i \kappa_i)$.  \hfill
\qed \medskip

We then have: 

\begin{prop} \label{free}
$(\gamma_{i_1}^{\al_1}\cdots \gamma_{i_p}^{\al_p})_{i_1<i_2<\cdots, 
\al_i \geq 1}$ is a free family of $\cA$. 
\end{prop}

{\em Proof.} There is an algebra map $i$ from $\cA$ to $FO$, mapping
each $\gamma[\eps]$ to $\eps$ in $FO^{(1)}$ (for $\eps = \sum_i \eps_i
z^i$, we set $\gamma[\eps] = \sum_i \eps_i \gamma_i$). This map also
sends the product $\gamma[\eps_1] \cdots \gamma[\eps_n]$ to
$S_q(\eps_1\otimes\cdots \otimes \eps_n)$. Let us show that the image
by $i$ of $(\g_{i_{1}}^{\al_{1}}\ldots
\g_{i_{p}}^{\al_{p}})_{i_{1}<\ldots< i_{p},\al_{i}\ge 1}$ forms a free
family in $FO$.

{}From Prop. \ref{str:FO} follows that multiplication by
$\prod_{i<j}q_-(z_i,z_j)$ -- call this map $m$ -- defines an
isomorphism of $FO^{(n)}$ with the image of the endomorphism $Sym_q$
of $k^{\hat\otimes n}[[\hbar]]$, defined by
$$Sym_q(\eps) = \sum_{\sigma\in
  \SS_n}\eps(z_{\sigma(1)},\cdots,z_{\sigma(n)}) \prod_{i<j, \sigma(i) >
  \sigma(j)}q_+(z_i,z_j)\prod_{i<j, \sigma(i) <
  \sigma(j)}q_-(z_i,z_j).$$ Moreover, we have $Sym_q(\eps_1 \otimes
\eps_2 \cdots) = m \circ i (\gamma[\eps_1]\gamma[\eps_2]\cdots)$.

Suppose now that some combination $\sum_{i_1\leq i_2 \leq \cdots}
\la_{i_1 i_2 \cdots} \gamma_{i_1}\gamma_{i_2} \cdots$ is zero in
$\cA$. It follows that the image by $Sym_q$ of the
combination $\sum_{i_1\leq i_2 \leq \cdots}\la_{i_1,i_2\cdots}\eps_{i_1}
\otimes \cdots$ is zero.  Let $\al$ be the smallest $\hbar$-adic
valuation of all $\la_{i_1,i_2,\cdots}$, then the leading term in $\hbar$
of this equality gives $\prod_{i<j}(z_i - z_j)\sum_{\sigma\in \SS_n}
\eps^{\sigma} = 0$, with $\eps = \sum_{i_1 \leq i_2 \cdots}\la_{i_1
  i_2\cdots}^{(\al)}\eps_{i_1} \otimes \eps_{i_2} \cdots$, which
implies the $\sum_{\sigma\in \SS_n} \eps^{\sigma} = 0$, so that
$\la_{i_1,i_2,\cdots}^{(\al)} = 0$, and all $\la_{i_1,i_2,\cdots}$ are $0$.
\hfill \qed \medskip

As a consequence of Lemma \ref{generating} and Prop. \ref{free}, we get 

\begin{prop} \label{PBW}
$(\gamma_{i_1}^{\al_1}\cdots \gamma_{i_p}^{\al_p})_{i_1<i_2<\cdots, 
\al_i \geq 1}$ is a topological basis of $\cA$. 
\end{prop}

We have also 

\begin{prop}
$i$ is an algebra isomorphism between $\cA$ and $FO$. 
\end{prop}

{\em Proof.} We have seen that $i(\gamma[\eps_1] \cdots
\gamma[\eps_n]) = S_q(\eps_1 \otimes \cdots \otimes \eps_n)$, so that
$i$ is surjective.  

{}From Prop. \ref{PBW} and the proof of Prop. \ref{free} follows that
the image by $i$ of a basis of $\cA$ is a free family of $FO$, so that
$i$ is injective.
\hfill \qed

\subsection{PBW result for $U_{\hbar}\G$}

Let $\cB^{+}$ be the quotient of the algebra $T(k)\hat{}[[\hbar]]$ by the
following relations: let $e'[\eps]$ denote the element of
$T(k)\hat{}[[\hbar]]$ 
corresponding to $\eps\in k$, and $e'(z)=\sum_{i\in
\ZZ}e'[\eps^{i}]\eps_{i}(z)$, the relations are the Fourier
coefficients of (\ref{e-e}), with the series replaced by their analogues
with primes. 

Similarly, let $\cB^{-}$ be the quotient of the algebra
$T(k)\hat{}[[\hbar]]$ by the
following relations: let $f'[\eps]$ denote the element of
$T(k){\hat{}}[[\hbar]]$ 
corresponding to $\eps\in k$, and $f'(z)=\sum_{i\in
\ZZ}f'[\eps^{i}]\eps_{i}(z)$, the relations are the Fourier
coefficients of (\ref{f-f}), with the series replaced by their analogues
with primes. 

Finally, let $\cB^{0}$ be the quotient of the algebra $T(k\oplus \CC K'
\oplus \CC D'){\hat{}}[[\hbar]]$ by the
following relations: let 
for $\eps\in k$, $h[\eps]$ denote the 
element of $T(k\oplus \CC K' \oplus \CC D'){\hat{}}[[\hbar]]$ 
corresponding to $\eps$,
the relations are (\ref{comm}),
(\ref{h-h}), (\ref{D-h}), with the $h^{+}[r], h^{-}[\la]$ replaced by 
$h'[r],h'[\la]$. 

There are algebra morphisms from 
$\cB^{\pm}$, $\cB^{0}$ to $U_{\hbar}{\G}$, associating to each
generator its version without prime.  

We then have 
\begin{lemma} \label{triangle}
The composition of the above algebra morphisms with the multiplication
of $U_{\hbar}{\G}$ defines a linear map
$$
i_{\G}:\cB^{+}\hat \otimes \cB^{0}\hat\otimes \cB^{-}\to U_{\hbar}{\G},
$$
which is a linear isomorphism (the tensor products are completed over
$\CC[[\hbar]]$).   
\end{lemma}

\noindent
{\em Proof.} We first consider the case of the algebra
$U_{\hbar}\G'$, defined as the algebra with the same generators (except
$D$) and relations as $U_{\hbar}\G$. Let $\cB^{\prime 0}$
be the analogue of algebra $\cB^{0}$ without generator $D'$. 
In that case, we obtain easily that 
$(b^{+}_{i}b^{\prime 0}_{j}b^{-}_{k})_{i,j,k}$ is a base of
$U_{\hbar}\G'$, if  
$(b^{+}_{i})_{i},(b^{\prime 0}_{j})_{j},(b^{-}_{k})_{k}$ are images of
bases of 
$\cB^{+},\cB^{\prime 0}, \bar \cB^{-}$.

Then we check that the r.h.s. of formulas (\ref{D-h}),
(\ref{D-ef}) define a derivation of $U_{\hbar}\G'$. We then apply the
PBW result for crossed products of algebras by derivations, and obtain
for $U_{\hbar}\G$ a base
$(b^{+}_{i}b^{\prime 0}_{j}b^{-}_{k}D^{s})_{i,j,k;s\ge 0}$. 
We finally
make use of (\ref{D-ef}), $x^{\pm}=f$, to pass (by triangular
transformations) from this base to
$(b^{+}_{i}b^{\prime 0}_{j}D^{s} b^{-}_{k})_{i,j,k;s\ge 0}$. 

Since $(b^{\prime 0}_{j}D^{s})_{j;s\ge 0}$ is the image of a base of
$\cB^{0}$, this final base has the desired form. 
\qed\medskip

\begin{lemma} \label{basis-for-B}
$\cB^{\pm}$ are topologically
spanned by 
$e'[\eps_{i_{1}}]^{\al_{1}} \cdots e'[\eps_{i_{p}}]^{\al_{p}}$, 
$i_{1}< \ldots <
i_{p}$, $\al_{i}\ge 1$, resp. 
$f'[\eps_{i_{1}}]^{\al_{1}} \cdots f'[\eps_{i_{p}}]^{\al_{p}}$, 
$i_{1}< \ldots <
i_{p}$, $\al_{i}\ge 1$.
\end{lemma}

\noindent{\em Proof.} This follows directly from the analogue of 
Prop. \ref{PBW} (where $\CC((\zeta))$ is replaced by a direct sum
$\oplus_{s\in S}\CC((\zeta_{s}))$), the
fact that $(z_{s}-w_{s})a_{0} \neq - (z_{s}-w_{s})a_{0}$, and the
computation of \cite{ER} preceding Thm. 5: 
\begin{align*}
  & e^{2\hbar[-\psi_{0}(\g,\pa_{z}\g,...)+\tau]}
  [z-w+\beta\psi_{-}(\g, \pa_{z}\g,...)] \cdot \\ &
  [e^{2\hbar[-\psi_{0}(\g,\pa_{z}\g,...)+\tau]} (z-w+\beta\psi_{-}(\g,
  \pa_{z}\g,...))]^{\widetilde{\ }} \\ & =
  [z-w+\beta\psi_{+}(\g,\pa_{z}\g,...)] \cdot 
  [z-w+\beta\psi_{+}(\g, \pa_{z}\g,...)]^{\widetilde{\ }}, 
\end{align*}
where $\beta(z,w) = (z-w)a_{0}$, 
and this is in turn written 
$$
e^{2\hbar(-\psi_{0}-\wt\psi_{0})}e^{2\hbar
{{T_{z}-T_{w}}\over{\pa_{z}+\pa_{w}}}(\g-\wt\g)}
{{1+ G \psi_{-}(\g,\pa_{z}\g,...)}\over{1+ G \psi_{+}(\g,\pa_{z}\g,...)}}
{{1- G 
\psi_{-}(\wt\g,\pa_{w}\wt\g,...)}\over{1- G \psi_{+}(\wt\g,
\pa_{w}\wt\g,...)}}=1,
$$
which amounts to the statement (3.17) of \cite{ER}. 
\qed\medskip

\begin{prop} \label{PBW-for-Ug}
The 
$$ 
e[\eps_{i_{1}}]^{\al_{1}} \cdots e[\eps_{i_{p}}]^{\al_{p}}
h[\eps_{k_{1}}]^{\gamma_{1}} \cdots h[\eps_{k_{r}}]^{\gamma_{r}}
D^{d} K^{t}
f[\eps_{j_{1}}]^{\beta_{1}} \cdots f[\eps_{j_{q}}]^{\beta_{q}}
$$
with $i_{1}< \ldots <
i_{p}$, $\al_{i}\ge 1$, 
$j_{1}< \ldots <
j_{q}$, $\beta_{i}\ge 1$, 
$k_{1}< \ldots <
k_{r}$, $\gamma_{i}\ge 1$, $d,t\ge 0$, 
form a topological basis of $U_{\hbar}\G$. 
\end{prop}

\noindent{\em Proof.}
This follows from Lemmas \ref{triangle}, \ref{basis-for-B} and the fact
that 
$h'[\eps_{k_{1}}]^{\gamma_{1}} \cdots h'[\eps_{k_{r}}]^{\gamma_{r}}
D^{d} K^{t}
$, $k_{1}< \ldots <
k_{r}$, $\gamma_{i}\ge 1$, $d,t\ge 0$, 
forms a basis of $\cB^{0}$. 
\qed\medskip

\begin{remark} It is straightforward to repeat the resoning above to
obtain topological bases of the $U_{\hbar}\G^{\hat\otimes N}$, as tensor
powers of the base obtained in Prop. \ref{PBW-for-Ug}. 
\qed
\end{remark}

\begin{remark} It would be interesting to explicitly compute, in the
case of the algebras $\cB^{\pm}$, the $\kappa_{0}$ provided by
Prop. \ref{PBW}.  
\end{remark}

\medskip\section{Subalgebra $U_{\hbar}\G_{R}$ of $U_{\hbar}\G$} \label{UgR}

\subsection{Presentation of $U_{\hbar}\G_{R}$} \label{UgR.1}

Recall that $\G$ contains as a Lie
subalgebra $\G_{R}$. In this section, we 
define a subalgebra $U_{\hbar}\G_{R}$ of 
$U_{\hbar}\G$, 
such that the inclusion 
$U_{\hbar}\G_{R}\subset U_{\hbar}\G$ is a deformation of 
$U\G_{R}\subset U\G$. 

Let $U_{\hbar}\G_{R}$ be the algebra with generators
$\wt D,\wt e[r],\wt f[r],\wt h[r],R,
r\in R$ and relations
\begin{equation} \label{lin-R}
\wt x[\al_{1}r_{1}+\al_{2}r_{2}]
=\al_{1}\wt x[r_{1}] +\al_{2}\wt x[r_{2}], \quad x=e,f,h,  
\end{equation}
\begin{equation} \label{h-hR}
[\wt h[r], \wt h[r']]=0, \quad \forall r,r'\in R,
\end{equation}
\begin{equation}\label{h-efR}
[\wt h[r], \wt e[r']]=2\wt e[rr'], \quad
[\wt h[r], \wt f[r']]=-2\wt f[rr'], \quad \forall r,r'\in R,
\end{equation}
\begin{align}  \label{e-eR}
\wt e[r_{1}\al]\wt e[r_{2}] & - \wt e[r_{1}]\wt e[\al r_{2}]
+ \wt e[r_{1}\psi_{-}^{(1)}\g(\al)^{(1)}]\wt
e[r_{2}\psi_{-}^{(2)}\g(\al)^{(2)}]
= 
\\ \nonumber 
&\wt e[r_{2}(q^{2(\tau-\phi)})^{(2)}]
\wt e[\al r_{1}(q^{2(\tau-\phi)})^{(1)}]
- 
\wt e[r_{2}(q^{2(\tau-\phi)})^{(2)} \al]
\wt e[r_{1}(q^{2(\tau-\phi)})^{(1)}]
\\ \nonumber 
&+
\wt e[r_{2}(\psi_{+}q^{2(\tau-\phi)})^{(2)}\gamma(\al)^{(2)}]
\wt e[r_{1}(\psi_{+}q^{2(\tau-\phi)})^{(1)}\gamma(\al)^{(1)}],
\end{align}

\begin{align} \label{f-fR}
\wt f[r_{2}]
\wt f[\al r_{1}]
& - 
\wt f[r_{2}\al]
\wt f[r_{1}]
+
\wt f[r_{2}(\psi_{-})^{(2)}\gamma(\al)^{(2)}]
\wt f[r_{1}(\psi_{-})^{(1)}\gamma(\al)^{(1)}]
=
\\ \nonumber 
& \wt f[(q^{2(\tau-\phi)})^{(1)}r_{1}\al]
\wt f[(q^{2(\tau-\phi)})^{(2)}r_{2}]  
- \wt f[(q^{2(\tau-\phi)})^{(1)}r_{1}]\wt f[(q^{2(\tau-\phi)})^{(2)}
\al r_{2}]
\\ \nonumber 
&
+ \wt f[r_{1}(q^{2(\tau-\phi)}\psi_{+})^{(1)}\g(\al)^{(1)}]\wt
f[r_{2}(q^{2(\tau-\phi)}\psi_{+})^{(2)}\g(\al)^{(2)}],
\end{align}
\begin{equation} \label{e-fR}
[\wt e[r_{1}],\wt f[r_{2}]]=\sum_{s\in
S}\res_{s}\{{1\over \hbar}
(r_{1}r_{2})(z)q^{((T+U)\wt h)(z)}\omega_{z}\}
\end{equation}
\begin{equation} \label{D-hR}
[\wt D, \wt h[r]]=\wt h[\pa r], 
\end{equation}
\begin{equation} \label{D-efR}
[\wt D, \wt x^{\pm}[r]] = \wt x^{\pm}[\pa r]+ {\hbar \over 2}
\sum_{i\in \NN} \wt h[e^{i}]\wt x^{\pm}[(A e_{i})r]
\end{equation}
for $x=e,f,h$; $x^{\pm}=e,f$; 
$\al,r_{i}\in R$, $\al_{i}$ scalars, $i=1,2$; $\wt
h(z)=\sum_{i\in \NN}\wt h[e^{i}]e_{i}(z)$; 
$$
\g(\al)=(\al\otimes 1 - 1 \otimes \al)a_{0}\in R\otimes R
$$ 
for $\al \in R$. 

\begin{notation} For $\xi\in R\otimes
R$, $\xi=\sum_{i}\xi_{i}\otimes \xi'_{i}$, and $a,b\in R$, we denote
$\sum_{i}\wt x[a\xi_{i}] \wt x[b\xi'_{i}]$
as
$\wt x[a\xi^{(1)}] \wt x[b\xi^{(2)}]$. The operator $A$ arising in
(\ref{D-efR}) has been defined in (\ref{A}). 
Note that the sums arising in (\ref{D-efR}) have only finite non-zero
terms, since $U$ has the property that for any sequence $(\xi_{i})$ with
$\xi_{i}\to 0$, $U\xi_{i}$ is zero for $i$ large enough, and both
sequences $e_{i},\pa e_{i}$ tend to zero; and on the other hand, $(\pa
e_{i})_{R}=0$ for $i$ large enough. 
\qed\medskip
\end{notation}

\begin{remark} {\it On relations (\ref{e-eR}) and (\ref{f-fR}).}
The
complicated-looking formulas (\ref{e-eR}) and (\ref{f-fR}) are simply
obtained by pairing the vertex relations (\ref{e-e-variant}) and
(\ref{f-f-variant}), for $\al\in R$, with an element of $R\otimes R$.
\qed\medskip
\end{remark}

\begin{remark} {\it Generating series for relations (\ref{e-eR}), 
(\ref{f-fR}) and (\ref{e-fR}).}
Let us introduce the generating series 
$$
\wt e(z)=\sum_{i\in \NN}\wt e [e^{i}]e_{i}(z), \quad 
\wt f(z)=\sum_{i\in \NN}\wt f [e^{i}]e_{i}(z).
$$
Relations (\ref{e-eR}) and
(\ref{f-fR}) can then be obtained as Fourier modes of 
\begin{align}  \label{gen-e-e-R}
\left( \right.
[ &  \al(z)-\al(w) +\psi_{-}(z,w)\gamma(\al)(z,w)]
\wt e(z)\wt e(w) 
\left. \right)_{\La,\La}
\\ \nonumber 
& =
\left(
q^{2(\tau-\phi)(z,w)}
[\al(z)-\al(w)+\psi_{+}(z,w)\gamma(\al)(z,w)]
\wt e(w)\wt e(z) 
\right)_{\La,\La}
, \quad \forall \al\in R, 
\end{align}
and
\begin{align}  \label{gen-f-f-R}
\left( \right.
q^{2(\tau-\phi)(z,w)}
[ & \al(z)-\al(w) +\psi_{+}(z,w)\gamma(\al)(z,w)]
\wt f(z)\wt f(w)
\left. \right)_{\La,\La}
\\ \nonumber 
&=
\left(
[\al(z)-\al(w) +\psi_{-}(z,w)\gamma(\al)(z,w)]
\wt f(w)\wt f(z) \right)_{\La,\La}, \quad \forall \al\in R, 
\end{align}
where the index $\La,\La$ has the following meaning: for any vector
space $V$ and $\xi \in V \otimes (k \bar\otimes k)$ (with $k\bar\otimes
k = \limm_{N} k/k_{N} \otimes k/k_{N}$),
$\xi_{\La,\La}=(id_{V}\otimes pr_{\La}\otimes pr_{\La})\xi$, where
$pr_{\La}$ denote the projection on the second summand of $k=R\oplus
\La$. 
In terms of these generating series, the relations (\ref{e-fR}) take the
form 
\begin{equation} \label{gen-e-f-R}
[\wt e(z),\wt f(w)]=\left( {1\over \hbar} q^{((T+U)\wt h)(z)}\delta(z,w)
\right)_{\La,\La},
\end{equation}
which can be rewritten as
$$
[\wt e(z),\wt f(w)]=\sum_{i\in \NN}\left( {1\over \hbar} q^{((T+U)\wt
h)(z)}e^{i}(z)\right)_{\La} e_{i}(w)
=
\sum_{i\in \NN}\left( {1\over \hbar} q^{((T+U)\wt
h)(w)}e^{i}(w)\right)_{\La} e_{i}(z)
$$
or in ``mixed'' form
$$
[\wt e(z),\wt f[r]]=[\wt e[r],\wt f(z)]=\left({1\over \hbar}
q^{((T+U)\wt h)(z)}r(z)\right)_{\La},
$$
for any $r\in R$. 
\qed\medskip
\end{remark}

\subsection{PBW result for $U_{\hbar}\G_{R}$ and inclusion in
$U_{\hbar}\G$} \label{UgR.2}

Let $\cB_{R}^{+}$ be the algebra
with generators $e''[r], r\in R$, and relations (\ref{lin-R}), with $\wt
x$ replaced by $e''$, and (\ref{e-eR}), with $\wt e$ replaced by $e''$;
let $\cB_{R}^{-}$ be the algebra
with generators $f''[r], r\in R$, and relations (\ref{lin-R}), with $\wt
x$ replaced by $f''$, and (\ref{f-fR}), with $\wt f$ replaced by $f''$;
and let $\cB_{R}^{0}$ be the algebra
with generators $D'',h''[r], r\in R$, and relations
(\ref{lin-R}), with $\wt x$ replaced by $h''$, (\ref{h-hR}) and
(\ref{D-hR}), 
with $\wt D,\wt h$ replaced by $D'',h''$.

\begin{lemma} \label{inject-BR-B}

1) There are injective algebra morphisms $i^{\pm}, i^{0}$ from 
$\cB_{R}^{\pm},\cB_{R}^{0}$ to $\cB^{\pm},\cB^{0}$,
sending each $x''[r]$ to $x'[r]$, and $D''$ to $D'$, $x=e,f,h,r\in
R$. 

2) Topological bases of $\cB_{R}^{\pm},\cB_{R}^{0}$ are given
by the $$(e''[e^{i_{1}}]^{\al_{1}}\ldots
e''[e^{i_{p}}]^{\al_{p}})_{i_{1}<\ldots < i_{p}, \al_{i}\ge 1},
\quad
(f''[e^{i_{1}}]^{\al_{1}}\ldots
f''[e^{i_{p}}]^{\al_{p}})_{i_{1}<\ldots < i_{p}, \al_{i}\ge 1},
$$ 
and
$(h''[e^{i_{1}}]^{\al_{1}}\ldots
h''[e^{i_{p}}]^{\al_{p}}D^{s})_{i_{1}<\ldots < i_{p}, \al_{i}\ge 1, s\ge 0}$.

\end{lemma}

\noindent{\em Proof.}
 1) (\ref{e-e-variant}) [resp. (\ref{f-f-variant})] with $e$ (resp. $f$)
replaced by $e'$ (resp. $f'$), are relations of $\cB^{+}$ (resp. of
$\cB^{-}$). Pair them with 
$r_{1}(z)r_{2}(w)$. We obtain relations (\ref{e-eR}), (\ref{f-fR}), with
$e',f'$ instead of $\wt e,\wt f$. This shows that the maps
$x''[r]\mapsto x'[r]$, $x=e,f$
extend to morphisms from $\cB_{R}^{\pm}$ to $\cB^{\pm}$. 
The statement on the map $h''[r]\mapsto h'[r], D''\mapsto D'$ is
evident.  

2) The case of $\cB_{R}^{0}$ is obvious. Let us treat $\cB_{R}^{+}$. 
From relations (\ref{e-eR}) follows
\begin{equation} \label{toto}
[e''[r_{0}\al],e''[r_{1}]]-[{e}''[r_{0}],{e}''[\al r_{1}]]
\in \hbar \cB_{R}^{+},
\end{equation}
for any $r_{0},r_{1},\al\in R$. 
We then get, setting $r_{0}=1$ in (\ref{toto}),  
$$
[{e}''[\al], {e}''[\beta]]\in 
[{e}''[1],{e}''[\al\beta]]+\hbar\cB_{R}^{+}; 
$$
on the other hand, setting $r_{0}=r_{1}=1$ in (\ref{toto}), we find that 
$$
2[{e}''[\al],{e}''[1]]\in \hbar\cB_{R}^{+}
$$
so that any commutator $[e''[\al],e''[\beta]]$ is in $\hbar\cB_{R}^{+}$. 
This shows that any monomial can be transformed into a combination of
the ${e}''[e_{i_{1}}]^{\al_{1}}\ldots
{e}''[e_{i_{p}}]^{\al_{p}}$, with $i_{1}<\ldots< i_{p}$, and
$\al_{i}\ge 1$.

Suppose now that some combination $\sum_{i\ge
0}\hbar^{i}\sum_{i_{1}<\ldots < i_{p}, \al_{i}\ge
1}a^{(i)}_{i_{j},\al_{j}}e''[e^{i_{1}}]^{\al_{1}}
\ldots e''[e^{i_{p}}]^{\al_{p}}$ is zero 
(where for each $i$, the second sum is a finite one). Applying $i^{+}$
to this identity, we obtain the identity in $\cB^{+}$ 
$$
\sum_{i\ge
0}\hbar^{i}\sum_{i_{1}<\ldots < i_{p}, \al_{i}\ge
1}a^{(i)}_{i_{j},\al_{j}}e'[e^{i_{1}}]^{\al_{1}}
\ldots e'[e^{i_{p}}]^{\al_{p}}=0,
$$
which implies that all $a^{(i)}_{i_{j},\al_{j}}$ are zero due to
Lemma \ref{basis-for-B}. This shows that 
$$
(e''[e^{i_{1}}]^{\al_{1}}\ldots
e''[e^{i_{p}}]^{\al_{p}})_{i_{1}<\ldots < i_{p}, \al_{i}\ge 1}
$$ 
is topologically free in 
$\cB_{R}^{+}$. Part 2 of the lemma follows for $\cB_{R}^{+}$. The case
of $\cB_{R}^{-}$ is similar. 
\qed\medskip

\begin{lemma} \label{triang-R}
There are injective algebra morphisms from $\cB_{R}^{\pm}$ and
$\cB_{R}^{0}$ to $U_{\hbar}\G_{R}$, sending each $x''[r]$ to $\wt x[r]$,
and $D''$ to $\wt D$, $x=e,f,h,r\in R$.
The composition of the tensor product of these morphisms, with the
multiplication of $U_{\hbar}\G_{R}$, induces a linear isomorphism
$$
i_{\G_{R}}: \cB_{R}^{+}\hat\otimes \cB_{R}^{0} \hat\otimes \cB_{R}^{-}
\to U_{\hbar}\G_{R}. 
$$
\end{lemma}

{\em Proof.} Let $U_{\hbar}\G_{R}'$ be the algebra with the same
generators (without $\wt D$) and relations as $U_{\hbar}\G_{R}$. 
We can prove by direct computation 
that the r.h.s. of relations (\ref{D-hR}), (\ref{D-efR})
define derivations of this algebra. 

The proof of the lemma is then identical to that of 
Lemma \ref{triangle}.  \qed\medskip

\begin{prop} The map 
sending each $\wt x[r]$ to $x[r]$, $x=e,f,h$, $r\in R$, and $\wt D$ to $D$,
extends to an injective algebra morphism from $U_{\hbar} \G_{R}$ to
$U_{\hbar}\G$.
\end{prop}

{\em Proof.} Let us first show that this map extends to an algebra
morphism for $U_{\hbar}\G_{R}$ to $U_{\hbar}\G$. 
For any $\al\in R$, (\ref{e-e-variant}) and
(\ref{f-f-variant}) are relations of $U_{\hbar}\G$. Pairing them with
$r_{1}(z)r_{2}(w)$, ($r_{1},r_{2}\in R$), we obtain relations (\ref{e-eR}),
(\ref{f-fR}), with $\wt e,\wt f$ replaced by $e,f$. Moreover,
(\ref{h-h}) is a relation of $U_{\hbar}\G$; pairing it with 
$r_{1}(z)r_{2}(w)$, we
obtain (\ref{h-hR}) with $\wt x$ replaced by $x$ for $x=e,f,h$. Finally,
pairing relations (\ref{D-ef}),  (\ref{D-h}) 
with $r(z)$, $r\in R$, we
obtain relations (\ref{D-efR}), (\ref{D-hR}), with $\wt x$ replaced by
$x$. This shows that the map $\wt x[r]\mapsto x[r], \wt D\mapsto D$
extends to an algebra morphism from $U_{\hbar}\G_{R}$ to $U_{\hbar}\G$,
that we will denote by $\iota$. 

We easily check that the diagram 

\medskip \noindent
\begin{equation} \label{diagram}
\setlength{\unitlength}{0.25mm}
\begin{array}{ccc}
\cB_{R}^{+}\otimes \cB_{R}^{0}\otimes \cB_{R}^{-} & \lrar{i_{\G_{R}}} &
U_{\hbar}\G_{R}  \\
\ldar{ i^{+} \otimes i^{0} \otimes i^{-}} &     & \ldar{\iota} \\
\cB^{+}\otimes \cB^{0}\otimes \cB^{-} & \lrar{i_{\G}} & U_{\hbar}\G
\end{array}
\end{equation}
commutes. 
By Lemmas (\ref{triangle}) and (\ref{triang-R}), the horizontal arrows
are vector spaces isomorphisms. From Lemma (\ref{inject-BR-B}) follows that
the left vertical arrow is injective. It follows that $\iota$ is also
injective.  
\qed\medskip

\subsection{Dependence of $U_{\hbar}\G_{R}$ in $\tau$ and $\La$}

In \cite{ER} we showed that the various algebras $U_{\hbar}\G$,
associated to different choices of $\La$ and $\tau$, are all
isomorphic. We are now going to show that the same is true for their
subalgebras $U_{\hbar}\G_{R}$. We will denote with a superscript 
$(\La,\tau)$ the objects associated with a choice $(\La,\tau)$. 

We will study two families of changes of the pair $(\La, \tau)$, that
will generate all possible changes. The first is to change $(\La,\tau)$
into $(\La,\tau')$, where $\upsilon = \tau'-\tau$ is an arbitrary
antisymmetric  
element in $R^{\otimes 2}[[\hbar]]$. To it is associated the map
$u:\La[[\hbar]] \to R[[\hbar]]$ defined by $u(\la) = \langle \upsilon, 1
\otimes \la \rangle_{k}$. The second family of changes is parametrized
by some continuous anti-self-adjoint linear map $r: \La\to R$ (or
equivalently, by some antisymmetric $r_{0}$ in $R^{\otimes 2}$). We then
define a new pair $(\bar\La, \bar \tau)$ by the formulas $\bar\La =
(1+r)\La$ and $\bar\tau = \tau -\sum_{i\in \NN} T(r(e_{i})) \otimes
e^{i}$. 

Recall now the results of \cite{ER}. 

\begin{prop} (see \cite{ER}) \label{depend}
1) There is an isomorphism
$i^{\tau,\tau'}$ from
$U_{\hbar}\G^{(\La,\tau')}$ to $U_{\hbar}\G^{(\La,\tau)}$, such that 
$$
i^{\tau,\tau'} (e^{(\La,\tau')} (z) ) = e^{{\hbar \over 2}(u
h^{+(\La,\tau)})(z)}e^{(\La,\tau)} (z), \quad
i^{\tau,\tau'} (e^{(\La,\tau')} (z) ) = 
f^{(\La,\tau)} (z)e^{{\hbar \over 2}(u h^{+(\La,\tau)})(z)},
$$
$$
i^{\tau,\tau'} (h^{+(\La,\tau')} (z) ) = h^{+(\La,\tau)} (z), \quad  
i^{\tau,\tau'}(D^{(\La,\tau')}) = D^{(\La,\tau)}. 
$$
2) There is an isomorphism $i^{\La,\bar\La}$ from
$U_{\hbar}\G^{(\La,\tau)}$ to $U_{\hbar}\G^{(\bar\La,\bar\tau)}$ such
that 
$$
i^{\La,\bar\La}(x^{(\La,\tau)}(z)) = x^{(\bar\La,\bar\tau)}(z), \quad
i^{\La,\bar\La}(D^{(\La,\tau)}) = D^{(\bar\La,\bar\tau)},
$$
$x=e,f,h^{+}$. 
\end{prop}

\begin{prop} Both maps $i^{\tau,\tau'}$ and $i^{\La,\bar\La}$
restrict to isomorphisms
of $U_{\hbar}\G^{(\La,\tau)}_{R}$ with $U_{\hbar}\G^{(\La,\tau')}_{R}$
and $U_{\hbar}\G^{(\bar\La,\bar\tau)}_{R}$. Therefore the algebras
$U_{\hbar}\G^{(\La,\tau)}_{R}$ are isomorphic for all choices of
$(\La,\tau)$. 
\end{prop}

{\em Proof.} For $r\in R$
\begin{align*}
&i^{\tau,\tau'} (e^{(\La,\tau')} [r]) =\sum_{s\in S}\res_{s} 
\left( r e^{{\hbar \over 2}(u
h^{+(\La,\tau)})}e^{(\La,\tau)}  \right)\omega \\ 
& =  
\sum_{s\in S}\res_{s} \sum_{n\ge 0}{1\over {n!}} \sum_{i_{1},
\ldots, i_{n} \in \NN} \left( {\hbar\over 2} \right)^{n}
h^{+(\La,\tau)}[e^{i_{1}}] \ldots h^{+(\La,\tau)}[e^{i_{n}}]
( r u(e_{i_{1}}) \ldots  u(e_{i_{n}} )
e^{(\La,\tau)}  )(z)\omega_{z}
\\ &
=\sum_{s\in S}\res_{s} \sum_{n\ge 0}{1\over {n!}} \sum_{i_{1},
\ldots, i_{n} \in \NN} \left( {\hbar\over 2} \right)^{n}
h^{+(\La,\tau)}[e^{i_{1}}] \ldots h^{+(\La,\tau)}[e^{i_{n}}]
e^{(\La,\tau)} [r\cdot
 u(e_{i_{1}} ) \ldots  u(e_{i_{n}} )] 
\end{align*}
so $i^{\tau,\tau'}(e^{(\La,\tau')}[r])\in
U_{\hbar}\G_{R}^{(\La,\tau)}$. 
The proof
that $i^{\tau,\tau'}(f^{(\La,\tau')}[r])\in
U_{\hbar}\G_{R}^{(\La,\tau)}$ is similar 
and in the case of $i^{\tau,\tau'}(h^{+(\La,\tau')}[r])$ the analogous
statement is obvious. 
The inverse of $i^{\tau,\tau'}$ is $i^{\tau',\tau}$; in particular, 
$i^{\tau,\tau'}$
is bijective. It is also an algebra morphism by
Prop. \ref{depend}. This shows 1). The proof of 2) is obvious. 
\qed\medskip

\begin{remark} There is another isomorphism $\bar i^{\tau,\tau'}$ from
$U_{\hbar}\G^{(\La,\tau)}$ to $U_{\hbar}\G^{(\La,\tau')}$, such that 
$$
\bar i^{\tau,\tau'} (e^{(\La,\tau')} (z) ) = e^{{\hbar \over 2}(u
h^{+})(z)}e^{(\La,\tau)} (z), \quad
\bar i^{\tau,\tau'} (e^{(\La,\tau')} (z) ) = 
\bar f^{(\La,\tau)} (z)e^{{\hbar \over 2}(u h^{+})(z)},
$$
$$
\bar i^{\tau,\tau'} (h^{+(\La,\tau')} (z) ) = h^{+(\La,\tau)} (z),
\quad 
\bar i^{\tau,\tau'}(D^{(\La,\tau')}) = D^{(\La,\tau)}. 
$$
It is easy to see that it also yields an isomorphism from
$U_{\hbar}\G^{(\La,\tau)}$ to $U_{\hbar}\G^{(\La,\tau')}$. 
\end{remark}

\subsection{$U_{\hbar}\G_{R}$ and $\Delta,\bar\Delta$}

Let us define $U_{\hbar}\G_{R}\hat\otimes U_{\hbar}\G$, resp. 
$U_{\hbar}\G\hat\otimes U_{\hbar}\G_{R}$, 
as the quotients of $T(\G_{R}\oplus\G)\hat{}[[\hbar]]$, resp. 
$T(\G\oplus\G_{R})\hat{}[[\hbar]]$ by the usual relations. These are
complete subalgebras of $U_{\hbar}\G\hat\otimes U_{\hbar}\G$. 

\begin{prop} \label{additional}
$$
\Delta(U_{\hbar}\G_{R}) \subset U_{\hbar}\G\hat\otimes U_{\hbar}\G_{R}, 
\quad 
\bar \Delta(U_{\hbar}\G_{R}) \subset U_{\hbar}\G_{R}\hat\otimes
U_{\hbar}\G.  
$$
\end{prop}

{\em Proof.} It is enough to check these statements for the generators
of $U_{\hbar}\G_{R}$. They are obvious for $D$ and $h^{+}[r]$. Moreover,
$$
\Delta(e[r])=\sum_{s\in S}\res_{s}\left 
(e(z)\otimes q^{((T+U)h^{+})(z)})r(z)\omega_{z}\right) + 1\otimes e[r]; 
$$
the first term of the r.h.s. of this equality can
be decomposed as a sum of terms, the second factors of which all lie in
the algebra generated by the $h^{+}[r'], r'\in R$. 

We also have 
\begin{align*}
\Delta(f[r]) & =f[r]\otimes 1 + \sum_{s\in S}\res_{s}\left( q^{-h^{-}(z)}
\otimes (q^{-K_{1}\pa}f)(z)\right)r(z)\omega_{z} \\
& = f[r]\otimes 1 + \sum_{s\in S}\res_{s} \left( 
\sum_{p\ge 0}\sum_{i_{1},\ldots, i_{p}\in \NN}{{(-\hbar)^{p}}\over{p!}}
h^{-}[e_{i_{1}}] \ldots h^{-}[e_{i_{p}}]
e^{i_{1}}(z) \ldots e^{i_{p}}(z)  \right. \\ & \left. 
\otimes (q^{-K_{1}\pa}f)(z)\right)r(z)\omega_{z} \\
& = f[r]\otimes 1 +
\sum_{p\ge 0}\sum_{i_{1},\ldots, i_{p}\in \NN}{{(-\hbar)^{p}}\over{p!}}
h^{-}[e_{i_{1}}] \ldots h^{-}[e_{i_{p}}]
\otimes 
f[q^{K_{1}\pa}(e^{i_{1}} \ldots e^{i_{p}}r)]. 
\end{align*}
All $f[\pa^{s}(e^{i_{1}} \ldots e^{i_{p}}r)]$ belong to
$U_{\hbar}\G_{R}$. This ends the proof of the proposition in the case of
$\Delta$. 

In the case of $\bar \Delta$, the proof is similar. 
\qed\medskip

\section{$F$, universal $R$-matrices and Hopf algebra pairings}

In what follows, we will denote by $U_{\hbar}\N_{\pm}$ the algebras
$\cB^{\pm}$. 

Recall that the Hopf algebras $(U_{\hbar}\G_+,\Delta)$ and
$(U_{\hbar}\G_-,\Delta')$, as well as $(U_{\hbar}\bar\G_+,\bar\Delta)$
and $(U_{\hbar}\bar\G_-,\bar\Delta')$, are dual. Denote by $\langle ,
\rangle_{U_{\hbar}\G}$ and $\langle , \rangle_{U_\hbar\bar\G}$ the
corresponding bilinear forms. They are defined by the formulas
\begin{equation} \label{elem}
\langle h^{+}[r], h^{-}[\la]\rangle_{U_{\hbar}\G} 
={2\over \hbar}\langle r, \la\rangle_{k}, \quad
\langle e[\eps], f[\eta]\rangle_{U_{\hbar}\G} 
={1\over \hbar}\langle \eps, \eta\rangle_{k}, 
\end{equation}
for $\eps,\eta\in k$, $r\in R,\la\in\La$,
$$ \langle D,K\rangle_{U_{\hbar}\G} = 1, \quad \langle D, \A(k)
\rangle_{U_{\hbar}\G}=\langle \A(k),K\rangle_{U_{\hbar}\G}=0,
$$
and 
\begin{equation} \label{bar:elem}
\langle h^{+}[r], h^{-}[\la] \rangle_{U_{\hbar}\bar{\G}} 
={2\over \hbar}\langle r, \la\rangle_{k}, \quad
\langle f[\eps], e[\eta]\rangle_{U_{\hbar}\bar{\G}} 
={1\over \hbar}\langle \eps, \eta\rangle_{k}, 
\end{equation}
for $\eps,\eta\in k$, $r\in R,\la\in\La$,
$$ 
\langle D,K\rangle_{U_{\hbar}\bar{\G}} = 1, \quad \langle D, \A(k)
\rangle_{U_{\hbar}\bar{\G}}=\langle
\A(k),K\rangle_{U_{\hbar}\bar{\G}}=0.
$$

{}From the Hopf algebra pairing rules follows immediately
\begin{lemma} 
  Let $U_{\hbar}\HH_{R}$ be the subalgebra of $U_\hbar\G$ generated by
  $D$ and the $h^+[r], r$ in $R$, and let $U_{\hbar}\HH_{\La}$ be the
  subalgebra of $U_\hbar\G$ generated by $K$ and the $h^+[\la], \la$ in
  $\La$.

  For any $x^{\pm}$ in $U_{\hbar}\N_{\pm}$ and $t_R,t_{\La}$ in
  $U_\hbar\HH_R$ and $U_{\hbar}\HH_{\La}$, we have
\begin{equation} \label{khor}
  \langle t_R x^+ , t_{\La}x^-\rangle_{U_\hbar\G} =
  \varepsilon(t_R)\varepsilon(t_\La) \langle x^+,
  x^-\rangle_{U_\hbar\G}
\end{equation}
and
\begin{equation} \label{pak}
  \langle x^-t_R, x^+t_{\La}\rangle_{U_\hbar\bar\G} =
  \varepsilon(t_R)\varepsilon(t_\La) \langle x^-,
  x^+\rangle_{U_\hbar\bar\G}.
\end{equation}
\end{lemma}

Then: 

\begin{prop} \label{prop:ident:pairing}
  The restrictions to $U_\hbar\N_+ \times U_\hbar\N_-$and $U_\hbar\N_-
  \times U_\hbar\N_+$ of the pairings $\langle , \rangle_{U_{\hbar}\G}
  $ and $\langle , \rangle_{U_\hbar\bar\G}$ coincide up to
  permutation.
\end{prop}

{\em Proof.} Fix $\eps_i,\eta_j$ in $k$, $i = 1, \ldots, n$, $j = 1,
\cdots, m$. Let us compute $$\langle \prod_{i=1}^n e[\eps_i],
\prod_{j=1}^m f[\eta_j] \rangle_{U_{\hbar}\G}$$ (we denote by
$\prod_{i\in I} x_i$ the product $x_{i_1}\cdots x_{i_p}$, where $I$ is a
set of integers $\{i_{\al}\}$ with $i_1<i_2 < \ldots$). It is clear that
this is zero if $n$ is not equal to $m$. Assume that $n = m$, then let
us compute the generating series $\langle \prod_{i=1}^n e[\eps_i],
\prod_{j=1}^m f(z_j)\rangle_{U_{\hbar}\G}$. By the Hopf algebra
pairing rules, it is equal to
\begin{align} \label{khar}
  & \sum_{\sigma\in \SS_n} \langle \otimes_{i=1}^n e[\eps_i] ,
  \otimes_{i=1}^n \lbrace \prod_{l=1,\sigma(l)<i}^{\sigma^{-1}(i)}
      q^{-h^-(z_l)}   \\ & \nonumber 
      f(z_{\sigma^{-1}(i)}) \prod_{l=\si^{-1}(i)+1,\sigma(l)<i}^{l}
      q^{-h^-(z_l)} \rbrace \rangle_{U_\hbar\G^{\otimes n}}.
\end{align} 
Set 
\begin{equation} \label{q}
q(z,w) = q^{2\sum_i ((T+U)e_i)(z)e^i(w)}.
\end{equation} 
We have $q^{h^-(w)} f(z) q^{-h^-(w)} = q(z,w)^{-1}f(z)$. Using this
equality and (\ref{khor}), we identify (\ref{khar}) with
\begin{align*} 
& \sum_{\si\in \SS_n} \langle \otimes_{i=1}^n e[\eps_i],
\otimes_{i=1}^n f(z_{\si^{-1}(i)}) \rangle_{U_\hbar\G^{\otimes n}}
\prod_{l<(i), \si(l)>\si(i)} q(z_i,z_l)^{-1}
\\ & = \sum_{\si\in \SS_n} \prod_{i=1}^n \eps_{\si(i)}(z_i)
\prod_{l<i , \si(l)>\si(i)} q(z_i,z_l)^{-1}; 
\end{align*}
each term of this sum belongs to $\CC((z_1))((z_2))\cdots((z_n))$. 
Therefore  
\begin{align} \label{explicit:pairing}
& \langle \prod_{i=1}^n e[\eps_i], \prod_{i=1}^n f[\eta_i]
\\ & \nonumber 
= \sum_{\sigma\in \SS_n} \res_{z_n\in S}\cdots \res_{z_1\in S}
 \sum_{\si\in \SS_n} \prod_{i=1}^n \eps_{\si(i)}(z_i)
\prod_{i=1}^n \eta_{i}(z_i)
\prod_{l<i, \si(l)>\si(i)} q(z_i,z_l)^{-1} 
\omega_{z_1}\cdots \omega_{z_n}, 
\end{align}
where $\res_{z\in S}\la_z$ means $\sum_{s\in S}\res_{s}\la_z$.  

In the case of $U_{\hbar}\bar\G$, one can lead the similar
computation for $$ \langle \prod_{i=1}^n f[\eta_i], \prod_{i=1}^n
e[\eps_i] \rangle_{U_\hbar\bar\G}.$$ In that case one uses the
identity $(K^+(z),f(w)) = q(z,w)^{-1}$, with $q(z,w)$ defined by
(\ref{q}); the result is the r.h.s. of (\ref{explicit:pairing}).
\hfill \qed \medskip

Let us define the completion $U_{\hbar}\G \bar\otimes U_{\hbar}\G$ as
follows. Let $I_{N}\subset U_{\hbar}\G$ be the left ideal generated by
the $x[\eps]$, $\eps\in \prod_{s\in S}z_{s}^{N}\CC[[z_{s}]]$.  Define
$U_{\hbar}\G \bar\otimes U_{\hbar}\G$ as the inverse limit of the
$U_{\hbar}\G^{\otimes 2} / I_{N} \otimes U_{\hbar} \G + U_{\hbar}\G
\otimes I_{N}$ (where the tensor products are $\hbar$-adically
completed). $U_{\hbar}\G \bar\otimes U_{\hbar}\G$ is clearly a
completion of $U_{\hbar}\G^{\hat\otimes 2}$. One defines similarly
$U_{\hbar}\G^{\bar\otimes n}$. 

\begin{defin} \label{def:def:F}
Let $(\al^i),(\al_i)$ be dual bases of $U_\hbar\N_+$ and $U_\hbar\N_-$. We
 set
\begin{equation} \label{def:F}
  F = \sum_i \al^i \otimes \al_i.
\end{equation}
\end{defin}

>From (\ref{explicit:pairing}) follows that $F$ belongs to 
$U_{\hbar}\G \bar\otimes U_{\hbar}\G$. 

By (\ref{elem}), we also have 
\begin{equation} \label{approx:F}
  F \in 1 + \hbar\sum_{i}e[\eps^i]\otimes f[\eps_i] + \sum_{j\geq 2}
U_{\hbar}\N_+^{[j]}\otimes U_{\hbar}\N_-^{[j]},
\end{equation}
where $U_{\hbar}\N_{\pm}^{[j]}$ are the degree $j$ homogeneous
components of $U_{\hbar}\N_{\pm}$ (where the $e[\eps]$ and $f[\eps]$
are given homogeneous degree $1$).

\begin{prop} \label{prop:R:mat}
  Set $q = e^\hbar$ and
$$ \cR = q^{D\otimes K} \exp\left( {\hbar\over 2} \sum_i h^+[e^i]
  \otimes h^-[e_i]\right) F,
$$ and
$$ \bar\cR = F^{(21)}q^{D\otimes K} \exp\left( {\hbar\over 2} \sum_i
  h^+[e^i] \otimes h^-[e_i]\right);
$$ then $\cR$ and $\bar\cR$ are the universal $R$-matrices of
$(U_{\hbar}\G,\Delta)$ and $(U_{\hbar}\G_-,\bar\Delta)$ (viewed as the
doubles of $(U_{\hbar}\G_+,\Delta)$and
$(U_{\hbar}\bar\G_+,\bar\Delta)$), respectively.
\end{prop}

{\em Proof.} These statements are equivalent to the following ones: 
$$
\langle id\otimes b_+, \cR \rangle_{U_{\hbar}\G} = b_+, \quad
\langle \cR, b_- \otimes id \rangle_{U_{\hbar}\G} = b_-, 
$$
for $b_\pm \in U_\hbar\G_\pm$, 
and 
$$ \langle id\otimes \bar b_+, \bar\cR \rangle_{U_{\hbar}\bar \G} =
\bar b_+, \quad \langle \bar \cR, \bar b_- \otimes id
\rangle_{U_{\hbar}\bar \G} = \bar b_-,
$$ for $\bar b_\pm \in U_\hbar\bar\G_\pm$. (Here and later, we will
set $\langle id \otimes a, b \otimes c \rangle = \langle a,c \rangle
b$, $\langle a \otimes id, b \otimes c \rangle = \langle a,b \rangle
c$, etc.)  Let us show the first statement.

Assume that $b_+$ has the form $t_R x_+$, with $t_R$ in
$U_{\hbar}\HH_R$ and $x_+$ in $U_{\hbar}\N_+$. Set $\cK = q^{D\otimes K}
\exp\left( {\hbar\over 2} \sum_i h^+[e^i] \otimes h^-[e_i]\right)$;
and set $\cK = \sum_{j}K_{j} \otimes K'_{j}$. Then 
$$
\langle id\otimes b_+, \cR \rangle_{U_{\hbar}\G} = 
\sum_{j,i} K_j \al^i \langle b_+ , K'_j\al_i \rangle_{U_{\hbar}\G}. 
$$ Now set $\Delta(b_+) = \sum b_+^{(1)}\otimes b_+^{(2)}$; we have
$\langle b_+ , K'_j\al_i \rangle_{U_{\hbar}\G} = \sum \langle
b_+^{(1)}, K'_j\rangle_{U_{\hbar}\G} \langle b_+^{(2)},
\al_i\rangle_{U_{\hbar}\G}$. Only the part of $\Delta(b_+)$ whose
first factors are of degree zero contribute to this sum. This part is
equal to $\Delta(t_R)(1\otimes x_+)$. Therefore, if $\Delta(t_R) =
\sum t_R^{(1)} \otimes t_R^{(2)}$, we have 
$$
\langle id\otimes b_+, \cR \rangle_{U_{\hbar}\G} = 
\sum_{j,i} K_j \al^i
\sum \langle
t_R^{(1)}, K'_j\rangle_{U_{\hbar}\G} \langle t_R^{(2)} x_+, 
\al_i\rangle_{U_{\hbar}\G}; 
$$
by (\ref{khor}), this is equal to 
$$
\sum_{j,i} K_j \al^i
\sum \langle
t_R , K'_j\rangle_{U_{\hbar}\G} \langle x_+, 
\al_i\rangle_{U_{\hbar}\G} 
= 
\sum_{j} K_j x_+
\langle t_R , K'_j\rangle_{U_{\hbar}\G}. 
$$ One easily checks that $\sum_{j} K_j\langle t_R ,
K'_j\rangle_{U_{\hbar}\G} = t_R$. Therefore, the last sum is equal to
$b_+$.

The proof of the other statements is similar. 
\hfill \qed \medskip 

\begin{lemma} \label{K:cocycle}
$\cK$ satisfies the cocycle identity 
$$
\cK^{12}(\bar\Delta \otimes 1)(\cK)
= 
\cK^{23}(1 \otimes \bar\Delta)(\cK). 
$$
\end{lemma}

{\em Proof.} $(\cB^0,\bar\Delta)$ is a Hopf subalgebra of
$(U_{\hbar}\G,\bar\Delta)$. It is easy to check that
$(U_{\hbar}\HH_R,\bar\Delta)$ and $(U_{\hbar}\HH_{\La},\bar\Delta')$
are dual Hopf algebra, and that the double of
$(U_{\hbar}\HH_R,\bar\Delta)$ is $(\cB^0,\bar\Delta)$. Moreover, $\cK$
represents the identity pairing between these algebras. The identities
of the Lemma are then consequences of the quasi-triangular identities.
\hfill \qed \medskip

\begin{lemma} \label{K:conj}
$\cK$ conjugates the coproducts $\bar\Delta$ and $\Delta'$, that is 
\begin{equation} \label{woll}
\Delta'(x) = \cK \bar\Delta(x) \cK^{-1}
\end{equation}
for any $x$ in $U_{\hbar}\G$. 
\end{lemma}

{\em Proof.} Let us first prove this identity for $x$ in $\cB^0$. The
$R$-matrix identity for $(\cB^0,\bar\Delta)$ says that $\bar\Delta'(x)
= \cK \bar\Delta(x) \cK^{-1}$ for $x$ in $\cB^0$. On the other hand,
the restrictions of $\Delta$ and $\bar\Delta$ to $\cB^0$ coincide.
This proves (\ref{woll}) in this case. 

Let us now treat the case where $x = e(z)$. 
Set $\cK_{0} = \exp({\hbar\over 2} \sum_i h^+[e^i]
\otimes h^-[e_i])$ and $\cK_D = q^{D\otimes K}$. Then $\cK =
\cK_D\cK_0$.

We have $\bar\Delta(e(z)) = (e\otimes
q^{-h^-})(q^{K_2\pa} z) + 1 \otimes e(z)$.  Then we have
$$
[\sum_i h^+[e^i] \otimes h^-[e_i], 1\otimes e(z)] = 
(T+U)(h^+(q^{K_2\pa}z))_{\La} \otimes e(z), 
$$ so that 
$$
\cK_0 (1\otimes e(z))\cK_0^{-1} = \exp(\hbar \sum_i
h^+[e^i]((T+U)(q^{K_2\pa}e_i))_{\La}(z)) \otimes e(z), 
$$
and 
\begin{align} \label{sznit}
  & \cK_D\cK_0 (1\otimes e(z))\cK_0^{-1}\cK_D^{-1} \nonumber \\ & =
  \nonumber \exp(\hbar \sum_i
  h^+[q^{K_2\pa}e^i]((T+U)(q^{K_2\pa}e_i))_{\La}(z)) \otimes e(z) \\ &
  = \nonumber \exp(\hbar \sum_i h^+[e^i](T+U)e_i(z)) \otimes e(z)\\ &
  = q^{(T+U)h^+}(z) \otimes e(z).
\end{align}

On the other hand, 
\begin{align} \label{wollis}
\cK_0 (e(z)\otimes 1)\cK_0^{-1} = e(z) \otimes q^{h^-(z)}, 
\end{align}
and 
$$ [h^-[e_i],h^-(z)] = {2\over \hbar} \left( ( q^{-K_2\pa}
  (T(q^{K_2\pa}e_i)_R) )(z) + Ue_i(z) - (q^{-K_2\pa} U((q^{K_2\pa}
  e_i)_{\La}) )(z) \right), 
$$
so that 
\begin{align*} & [{\hbar\over 2}\sum_i h^+[e^i] \otimes h^-[e_i], 
  1\otimes h^-(z)] \\ & = 
  (q^{-K_2\pa} T (q^{K_2\pa}h^+)_R)(z) + Uh^+(z) -
  q^{-K_2\pa}(U(q^{K_2\pa}h^+)_{\La})(z) \otimes 1, 
\end{align*}
so that 
\begin{align*} & \cK_0 (1\otimes h^-(q^{K_2\pa}z)) \cK_0^{-1} 
\\ & = 1 \otimes
q^{K_2\pa}h^-(z) + \left( T(q^{K_2\pa}h^+)_R(z) + Uh^+(q^{K_2\pa}z) -
  (U(q^{K_2\pa}h^+)_{\La})(z) \right) \otimes 1
\end{align*}
and 
\begin{align} \label{hofen}
& \cK_0(1\otimes q^{-h^-(q^{K_2\pa} z)})\cK_0^{-1} \\ & = \nonumber  
\exp(-\hbar \left( T(q^{K_2\pa}h^+)_R(z) + Uh^+(q^{K_2\pa}z) -
  (U(q^{K_2\pa}h^+)_{\La})(z) \right) ) \otimes 
q^{-h^-(q^{K_2\pa} z)}. 
\end{align}
Finally, the product of (\ref{wollis}) and (\ref{hofen}) gives
\begin{align} \label{wollishofen}
& \cK_0 ((e\otimes q^{-h^-})(q^{K_2\pa}z)) \cK_0^{-1} \\ & = \nonumber 
\exp(-\hbar \left( T(q^{K_2\pa}h^+)_R(z) + Uh^+(q^{K_2\pa}z) -
  (U(q^{K_2\pa}h^+)_{\La})(z) \right) ) e(q^{K_2\pa}z)
\otimes 1. 
\end{align}
On the other hand, we have 
$$
(D + \pa_z - Ah^+(z)) e(z) = e(z) (D + \pa_z),  
$$
(an identity of differential operators). Therefore we have  
\begin{align} \label{lea}
q^{K_2(D + \pa_z - Ah^+(z))} e(z) = e(z) q^{K_2(D + \pa_z)},  
\end{align}

We are now in position to apply the following identity:
\begin{lemma} \label{tawil}
Let $L,f$ be elements of a Lie algebra, such that all $\ad^{k}(L)(f)$
commute with $f$. We have the identity, in the associated formal group  
$$
\al^{L+f} = \al^{L} \exp\left( {{1-\al^{-\ad L}} \over {\ad L}}(f) \right)
, 
$$
where $\al = e^{\hbar k}$, $k$ scalar and $\hbar$ a formal
parameter. 
\end{lemma}

{\em Proof. } Enlarge the Lie algebra by adjoining to it an element $F$, such
that $[L,F] = f$. We have in the associated formal group, 
\begin{align*}
\al^{L+f}  & = \al^{\Ad(e^{-F})(L)} = \Ad(e^{-F})(\al^{L}) = 
\al^{ L }
[ (\al^{-\ad(L)} ( e^{-F} ) )
e^{F} ] 
 \\  & = \al^{L} \exp(( 1 - \al^{ - \ad(L)})(F))  = \al^{L}  
\exp \left( { { 1 - \al^{-\ad(L)} } \over {\ad(L)} }(f)  
 \right) .  \end{align*}
\qed \medskip

Set in (\ref{lea}), $\al = q^{K_2}$, $L = D + \pa_z$ and $f =
-Ah^+(z)$. Then we need to compute 
$$
{{ 1 - q^{ - K_2\ad(D+ \pa_z)} }\over{\ad(D+\pa_z)}}(-A h^+)(z). 
$$
For this, we show: 
\begin{lemma} \label{haim}
Let $B$ be some linear operator from $\La[[\hbar]]$ to $R[[\hbar]]$,
such that $Be_{i}\to 0$ when $i\to\infty$, then

1) $\ad(D+ \pa_z)(Bh^{+}(z) ) 
= Ch^{+}(z)$, 
where $C = \pa \circ B - B \circ pr_{\La} \circ \pa$, 

2) $\Ad(\al^{D + \pa_z})( Bh^{+}(z)) = Gh^+(z)$, where
$G=\al^{\pa}\circ B \circ pr_{\La}\circ \al^{-\pa}$, and $\al =
e^{\hbar k}$, $k$ scalar.

3) assume that $B = \pa \circ E - E \circ pr_{\La} \circ \pa$, with $E$
a linear operator from $\La[[\hbar]]$ to $k[[\hbar]]$, then 
$$
{{\al^{\ad (D+\pa_z)}-1}\over{\ad (D+\pa_z) } }  
\left( Bh^{+}(z) \right) = 
Fh^{+}(z), 
$$
where $F = \al^{\pa}\circ E \circ pr_{\La}\circ \al^{\pa} - E$. 
\end{lemma}

{\em Proof.} 1) is obvious. 2) is obtained by first computing
$$
\ad^{k}(D+\pa_z)( Bh^{+}(z)[e^{i}]), 
$$ 
using 1). 
We again use this expression to obtain 3). 
\qed\medskip

Applying Lemma \ref{haim} with $B = A$, $E = \hbar(T+U)$, $\al = q^{-K_2}$,
we get
$$
{{ 1 - q^{ - K_2\ad(D+ \pa_z)} }\over{\ad(D+\pa_z)}}(-A h^+)(z). 
= {{ q^{ - K_2\ad(D+ \pa_z)} - 1}\over{\ad(D+\pa_z)}}(A h^+)(z). 
= Fh^+(z), 
$$
with $F = \hbar(q^{-K_2\pa}\circ (T+U) \circ pr_{\La} - (T+U))$. 
Therefore Lemma \ref{tawil} gives 
\begin{align*} 
& q^{K_2(D+\pa_z - Ah^+(z))} = q^{K_2(D+\pa_z} \exp( {{ 1 - q^{ -
      K_2\ad(D+ \pa_z)} }\over{\ad(D+\pa_z)}}(-A h^+)(z) )
\\ &  =
q^{K_2(D+\pa_z} \exp([q^{-K_2\pa}\circ \hbar(T+U) \circ pr_{\La} -
\hbar(T+U)]h^+(z)). 
\end{align*}

Finally, (\ref{lea}) gives 
$$
q^{K_2(D+\pa_z)} \exp([q^{-K_2\pa}\circ \hbar(T+U) \circ pr_{\La} -
\hbar(T+U)]h^+(z))e(z) = e(z) q^{K_2(D + \pa_z)},
$$
so that 
$$
q^{-K_2 D}e(z) q^{K_2D} = 
q^{K_2\pa_z} \{ \exp([q^{-K_2\pa}\circ \hbar(T+U) \circ pr_{\La} -
\hbar(T+U)]h^+(z))e(z)\}, 
$$
which coincides with the r.h.s. of (\ref{wollishofen}). Therefore 
\begin{equation} \label{man}
\cK_D\cK_0 ((e\otimes q^{-h^-})(q^{K_2\pa}z)) (\cK_D\cK_0)^{-1} 
= e(z) \otimes 1. 
\end{equation}

Adding up (\ref{sznit}) and (\ref{man}), we get 
$$
\Delta'(e(z)) = \cK \bar\Delta(e(z)) \cK^{-1}.  
$$

The proof is similar when $e(z)$ is replaced by the case of $f(z)$. 
\hfill \qed \medskip

\begin{prop} \label{prop:F:cocycle}
$F$ satisfies the cocycle identity
\begin{equation} \label{cocycle:id}
F^{12}(\Delta\otimes 1)(F) = F^{23}(1 \otimes \Delta)(F), 
\end{equation} 
\end{prop}

{\em Proof.} The quasi-triangularity identities for $\bar\cR$ imply that
$$ \bar\cR^{12}(\bar\Delta \otimes 1)(\bar\cR) = \bar\cR^{23}(1
\otimes \bar\Delta)(\bar\cR).
$$
Therefore, we have 
$$ F^{21}\cK^{12}(\bar\Delta\otimes 1)(F^{21})(\bar\Delta \otimes 1)(K)
= F^{32}\cK^{23}(1\otimes \bar\Delta)(F^{21})(1 \otimes \bar\Delta)(K).
$$
By Lemma \ref{K:conj}, it follows that 
$$ 
F^{21}(\Delta'\otimes 1)(F^{21})\cK^{12}(\bar\Delta \otimes 1)(\cK)
= F^{32}(1\otimes \Delta')(F^{21})\cK^{23}(1 \otimes \bar\Delta)(\cK).
$$
Lemma \ref{K:cocycle} then implies that 
$$
F^{21}(\Delta'\otimes 1)(F^{21})
= F^{32}(1\otimes \Delta')(F^{21}), 
$$ 
which is a the same as (\ref{cocycle:id}), up to permutation of
factors.  \hfill \qed \medskip

\begin{prop} \label{prop:F:twist}
  $\Delta$ and $\bar\Delta$ are conjugated by $F$: we have
  $\bar\Delta(z) = F \Delta(z) F^{-1}$, for any $x$ in $U_{\hbar}\G$.
\end{prop}

{\em Proof.} We know that $\Delta' = \cR \Delta \cR^{-1}$. Since $\cR
= \cK F$, it follows that $F\Delta F^{-1}= \cK^{-1} \Delta' \cK$. But by
Lemma \ref{K:conj}, $\cK^{-1} \Delta' \cK$ coincides with $\bar\Delta$. 
\hfill \qed \medskip

\begin{remark}
It is a general principle in $R$-matrix computations (see \cite{KT})
that factors of the $R$-matrix are also twists relating quantizations
of conjugated Manin triples. We see that this principle also holds
in our situation. \qed\medskip 
\end{remark}

\subsection{Orthogonals of $\cB_{R}^{\pm}$}

\begin{prop} \label{orthBR}
1)
The orthogonal of $\cB_{R}^{-}$ in $\cB^{+}$ for 
$\langle, \rangle_{U_{\hbar}\G}$ is the span $\N_{+}(R)\cB^{+}$
of all $e[r]e[\eps_{1}]\ldots
e[\eps_{p}]$, $p\ge 0$, $\eps_{i}\in k$, $r\in R$; 

2)
the orthogonal of $\cB_{R}^{+}$ in $\cB^{-}$ for 
$\langle, \rangle_{U_{\hbar}\G}$ is the span $\cB^{-} \N_{-}(R)$
of all 
$$
f[\eta_{1}]\ldots
f[\eta_{p}]f[r] , \quad p\ge 0, \eta_{i}\in k, r\in R.
$$
\end{prop}

{\em Proof.}
Let us compute $\langle e[r]e[\eps_{1}]\ldots
e[\eps_{p}], f[r_{1}]\ldots f[r_{p+1}] \rangle_{U_{\hbar}\G}$, 
$\eps_{i}\in k$, $r, r_{i}\in R$. 
Expand it as 
$$
\langle e[r] \otimes 
e[\eps_{1}]\ldots
e[\eps_{p}], \Delta' ( f[r_{1}]\ldots f[r_{p+1}] )
\rangle_{U_{\hbar}\G^{\hat\otimes 2}}
$$ 
(recall that $\Delta'$ is $\Delta$ composed with the exchange of
factors). According to the proof of Prop. \ref{additional},  
$\Delta' ( f[r_{1}]\ldots f[r_{p+1}] )$ can be decomposed as a sum of
terms, the first factors all of which are of the form $f[r'_{1}]\ldots
f[r'_{s}]$, $r'_{i}\in R$. Since $\langle e[r],f[r'_{1}]\ldots
f[r'_{s}]\rangle_{U_{\hbar}\G}$ is always zero (either because $s\ne 1$
of by isotropy of $R$), 
$$
\langle e[r]e[\eps_{1}]\ldots
e[\eps_{p}], f[r_{1}]\ldots f[r_{p+1}] \rangle_{U_{\hbar}\G}=0.
$$ 
This shows that $\N_{+}(R)\cB^{+} \subset (\cB_{R}^{-})^{\perp}$. 

To show that $\N_{+}(R)\cB^{+}$ is the whole of 
$(\cB_{R}^{-})^{\perp}$, let us consider the classical limit of the
situation. $\cB_{R}^{-}\subset \cB^{-}$ is a flat deformation of the
inclusion of symmetric algebras $S^{*}(R) \subset S^{*}(k)$. On the
other hand, the inclusion $\N_{+}(R)\cB^{+}\subset \cB^{+}$ is a flat
deformation of $S^{*}(k)R \subset S^{*}(k)$. Finally, let 
$\cB_{p}^{i}$ be 
the completion of the span of all $x[\eps_{1}]\ldots x[\eps_{p}],
\eps_{i}\in k$, $x=e,f$ for $i=+,-$,  and  multiply
$\langle, \rangle_{U_{\hbar}\G}$ 
by $\hbar^{p}$ on $\cB_{p}^{+}\times \cB_{q}^{-}$. Then the resulting
pairing is a deformation of the direct sum of the symmetric powers of
the pairing $\langle,
\rangle_{k}$. Since the orthogonal of $S^{*}(R)$ for this pairing is
$S^{*}(k)R$ (because $R$ is maximal isotropic, see Lemma
\ref{statement-R}), the orthogonal of $\cB_{R}^{-}$ cannot be larger than 
$\N_{+}(R)\cB^{+}$. This finally shows 1). 

Let us pass to the proof of 2). Let us compute 
$$
\langle e[r_{1}]\ldots
e[r_{p+1}], f[\eta_{1}]\ldots f[\eta_{p}] f[r] \rangle_{U_{\hbar}\G},
\quad r,r_{i}\in R, \eta_{i}\in k.
$$ 
Expand it as
$\langle \Delta ( e[r_{1}]\ldots
e[r_{p+1}] ) , f[\eta_{1}]\ldots f[\eta_{p}] \otimes  f[r]
\rangle_{U_{\hbar}\G^{\hat\otimes 2}}$. From the proof of
Prop. \ref{additional} follows that 
 $\Delta ( e[r_{1}]\ldots
e[r_{p+1}] )$ can be decomposed as a sum of terms, the second factors
all of which lie in the algebra generated by the $e[r], h^{+}[r'],
r,r'\in R$.  

Further decompose each of these second factors
as a sum of terms of the form
$$
e[r''_{1}]\ldots e[r''_{s}] h^{+}[\bar r_{1}]\ldots h^{+}[\bar r_{t}],
\quad r''_{i},\bar r_{i}\in R.
$$ 
The pairing of this with $f[r]$ gives
$$
\langle e[r''_{1}]\ldots e[r''_{s}]
\otimes h^{+}[\bar r_{1}]\ldots h^{+}[\bar r_{t}], \Delta'(f[r]) 
\rangle_{U_{\hbar}\G^{\hat\otimes 2}}.
$$ 
$\Delta'(f[r])$ is equal to the sum of $1\otimes f[r]$ and of 
a sum of terms, the first factors
of which are either $1$ or of the form $f[\rho], \rho\in R$. 

Since $\langle e[r''_{1}]\ldots e[r''_{s}], \rangle
f[\rho]\rangle_{U_{\hbar}\G}=0$ (either by degree reasons if $s\ne 1$ or
by isotropy of $R$), the only possibly non-trivial contribution is that
of 
$$
\langle e[r''_{1}]\ldots e[r''_{s}]
\otimes h^{+}[\bar r_{1}]\ldots h^{+}[\bar r_{t}], 1 \otimes f[r]
\rangle_{U_{\hbar}\G^{\hat\otimes 2}};
$$
but the pairing of $f[r]$ with
any $h^{+}[\bar r_{1}]\ldots h^{+}[\bar r_{t}]$ is zero.

This shows that  $\langle e[r_{1}]\ldots
e[r_{p+1}], f[\eta_{1}]\ldots f[\eta_{p}] f[r]\rangle_{U_{\hbar}\G}=0$, 
for any $r,r_{i}\in R, \eta_{i}\in k$, so that  $\cB^{-}\N_{-}(R)  
\subset (\cB_{R}^{+})^{\perp}$.

The proof that $\cB^{-}\N_{-}(R)$ is actually equal to 
$(\cB_{R}^{+})^{\perp}$ is similar to the argument used in the proof of
1). 
\qed
\medskip

\medskip\section{Quasi-Hopf structures}

\subsection{Factorization of $F$} \label{fctr}

We now recall our aim. We would like to decompose $F$ defined in
(\ref{def:F}) as a product
\begin{equation} \label{decomp}
F_{2}F_{1},\quad \on{ with} \quad  F_{1}\in U_{\hbar}\G\hat\otimes
U_{\hbar}\G_{R}, \quad 
F_{2}\in U_{\hbar}\G_{R}\hat\otimes U_{\hbar}\G.
\end{equation}
The interest of this decomposition lies in the following proposition.  

\begin{prop} \label{why-decomp}
For any decomposition (\ref{decomp}), the map $\Ad
(F_{1})\circ\Delta$ defines an algebra morphism from $U_{\hbar}\G_{R}$ to
$U_{\hbar}\G_{R}\hat\otimes U_{\hbar}\G_{R}$ (where the tensor product
is completed over $\CC[[\hbar]]$). 
\end{prop}

{\em Proof.}
Since $\bar \Delta=\Ad(F)\circ \Delta$, we have 
$$
\Ad (F_{1})\circ\Delta = \Ad (F_{2}^{-1}) \circ \bar\Delta. 
$$
The first map sends $U_{\hbar}\G_{R}$ to $U_{\hbar}\G\hat\otimes
U_{\hbar}\G_{R}$, 
and the second one to $U_{\hbar}\G_{R}\hat\otimes
U_{\hbar}\G$, by Prop. (\ref{additional}) and (\ref{decomp}); so both
maps send $U_{\hbar}\G_{R}$ to 
$ ( U_{\hbar}\G\hat\otimes
U_{\hbar}\G_{R} )  \cap  ( U_{\hbar}\G_{R}\hat\otimes
U_{\hbar}\G ) =U_{\hbar}\G_{R}\hat\otimes U_{\hbar}\G_{R}$. 
\qed\medskip

Let us now try to decompose $F$ according to (\ref{decomp}). Let
$(m_{i})$, resp. $(m'_{i})$ be a basis of $U_{\hbar}\G$ as a left,
resp. right $U_{\hbar}\G_{R}$-module. Assume $m_{0}=m'_{0}=1$. 
Due to the form of $F_{1}$ and $F_{2}$, we have decompositions 
$$
F_{2}=\sum_{i}(1 \otimes m'_{j}) F^{(j)}_{2}, \quad
F_{1}=\sum_{i}F^{(i)}_{1}(m_{i} \otimes 1), \quad 
F^{(i)}_{1} , F^{(j)}_{2}\in U_{\hbar}\G_{R}^{\hat\otimes 2}.  
$$
It follows that we have 
\begin{equation} \label{key}
F=\sum_{i}F_{2}F_{1}^{(i)}(m_{i} \otimes 1)
=
\sum_{j}( 1 \otimes m'_{j})F_{2}^{(j)}F_{1}. 
\end{equation}
Let now $\Pi$, resp. $\Pi'$ be the left, resp. right
$U_{\hbar}\G_{R}$-module morphisms from $U_{\hbar}\G$ to
$U_{\hbar}\G_{R}$, such that $\Pi(m_{i})=0$ for $i\neq 0$, $\Pi(1)=1$,
and $\Pi'(m'_{i})=0$ for $i\neq 0$, $\Pi'(1)=1$. 

>From (\ref{key}) follows that we should have 
\begin{equation} \label{possible}
F_{2}F_{1}^{(0)}=(\Pi\otimes 1)F, \quad 
F_{2}^{(0)}F_{1}=(1 \otimes \Pi')F. 
\end{equation}

We may and will assume that $m_{i}$, resp. $m'_{i}$ contains a basis of
$\cB^{+}$ as a left $\cB_{R}^{+}$-module, resp. of 
$\cB^{-}$ as a right $\cB_{R}^{-}$-module. Then, $\Pi$ maps $\cB^{+}$ to
$\cB_{R}^{+}$, and $\Pi'$ maps $\cB^{-}$ to $\cB_{R}^{-}$. It follows
that $[(\Pi\otimes 1)F]^{-1} F [(1 \otimes \Pi')F]^{-1}$ belongs to
$\cB^{+}\hat\otimes \cB^{-}$. 

Equation (\ref{possible}) determines the possible values of 
$F_{1}$ and $F_{2}$, up to 
right, resp. left multiplication by elements of
$U_{\hbar}\G_{R}^{\hat\otimes 2}$. 

\begin{prop} Let $F_{\Pi,\Pi'} = [(\Pi\otimes 1)F]^{-1} F [(1 \otimes
\Pi')F]^{-1}$; then 
\begin{equation} \label{main}
F_{\Pi,\Pi'} \in U_{\hbar}\G_{R}^{\hat\otimes 2}. 
\end{equation}
\end{prop}

{\em Proof.} Since $(\Pi\otimes 1)F\in U_{\hbar}\G_{R}\hat \otimes
U_{\hbar}\G$, and $(1\otimes \Pi')F\in U_{\hbar}\G\hat \otimes
U_{\hbar}\G_{R}$, 
(\ref{main}) is equivalent to showing that
\begin{equation}  \label{equiv}
F^{-1}[(\Pi\otimes 1)F]\in \cB^{+}\hat\otimes \cB^{-}_{R}, \quad
[(1 \otimes \Pi')F]F^{-1}\in \cB^{+}_{R}\hat\otimes \cB^{-}.  
\end{equation}

By Prop. \ref{prop:R:mat}, the universal $R$-matrices of $U_{\hbar}\G$
and $U_{\hbar}\bar{\G}$, $\cR$ and $\bar \cR$ are such that $\cR \in
(U_{\hbar}\G_{R}\hat\otimes U_{\hbar}\G)F$, and $\bar \cR^{21}\in
F(U_{\hbar}\G \hat\otimes U_{\hbar}\G_{R})$.

Since $\Pi$, resp. $\Pi'$ is a left, resp. right
$U_{\hbar}\G_{R}$-module morphism, it follows that  
\begin{equation} \label{F=R}
F^{-1}[(\Pi\otimes 1)F] = \cR^{-1}[(\Pi\otimes 1)\cR], \quad \on{and} \quad
[(1\otimes \Pi')F]F^{-1} = [(1\otimes \Pi')\bar \cR^{21}](\bar
\cR^{21})^{-1}.
\end{equation}

We will need the following result. 
\begin{lemma} \label{R=skew}
1) 
Let $S'$ denote the skew antipode of $U_{\hbar}\G$, then for any $x\in
\cB^{+}$, 
$$
\langle \cR^{-1}, id \otimes x \rangle_{U_{\hbar}\G} = S'(x). 
$$
2) Recall $\bar S$ is the antipode of $U_{\hbar} \bar{\G}$. For any
$y\in \cB_{-}$, we have
$$
\langle \bar \cR^{-1}, id \otimes y \rangle_{U_{\hbar}\bar{\G}} = \bar S(y). 
$$
\end{lemma}

{\em Proof.} Since $(\Delta\otimes 1)(\cR^{-1})=
(\cR^{-1})^{23} (\cR^{-1})^{13}$, $\sigma': \cB_{+} \to
U_{\hbar}\G_{+}$, $x\mapsto \langle
\cR^{-1}, id \otimes x \rangle_{U_{\hbar}\G}$ is an algebra morphism from
$\cB^{+}$ to $U_{\hbar}\G_{+}^{opp}$ (that is, $U_{\hbar}\G_{+}$ with
the opposite multiplication). Since $S'$ is also an algebra morphism from
$\cB^{+}$ to $U_{\hbar}\G_{+}^{opp}$, it suffices to check that $S'$ and
$\sigma'$ coincide on the generators of $\cB^{+}$.  

>From (\ref{Delta-e}) follows that $S'(e(z)) = -
e(z)q^{-((T+U)h^{+})(z)}$, on the other hand
\begin{align*}
\sigma'(e (z) )  & = \langle \cR^{-1}, id \otimes e (z) 
\rangle_{U_{\hbar}\G} 
\\ &
=
\langle ( 1 -\hbar \sum_{i\in\ZZ} e[\eps^{i}]\otimes f[\eps_{i}] + \cdots) 
e^{-{\hbar \over 2} \sum_{j\in\NN} h^{+}[e^{j}]\otimes h^{-}[e_{j}]}
, id 
\otimes e (z)
\rangle_{U_{\hbar}\G} 
\\ & 
= 
\langle (-\hbar \sum_{i\in\ZZ} e[\eps^{i}]\otimes f[\eps_{i}] ) 
e^{-{\hbar \over 2} \sum_{j\in\NN} h^{+}[e^{j}]\otimes h^{-}[e_{j}]}
, id 
\otimes e (z) 
\rangle_{U_{\hbar}\G} 
\\ &
=
\langle (-\hbar) \sum_{i\in\ZZ} e[\eps^{i}]^{[1]} f[\eps_{i}]^{[2]}
e^{-{\hbar \over 2} \sum_{j\in\NN} h^{+}[e^{j}]^{[1]} h^{-}[e_{j}]^{[3]}}
, id 
\otimes e (z) \otimes q^{((T+U)h^{+})(z)} 
\rangle_{U_{\hbar}\G} 
\\ & 
= - e(z)q^{-((T+U)h^{+})(z)}; 
\end{align*}
in the r.h.s. of the first equality, the 
notation means that we perform the pairing of the second factors of
a decomposition of $\cR^{-1}$ with $e(z)$, so this r.h.s.
belongs to $U_{\hbar}\G_{R}$; 
the third equality is because $e(z)$ cannot have a nontrivial pairing
with a product of zero or more than two $f[\eps_{i}]$'s; in the fourth
equality, we used the notation $a^{[i]}$ for 
$1^{\otimes (i-1)} \otimes a \otimes 1^{\otimes (3-i)}$; the last
equality follows from 
$$
\langle (-\hbar)\sum_{i\in\ZZ}e[\eps^{i}]
\otimes f[\eps_{i}], id \otimes  e(z)\rangle_{U_{\hbar}\G} = - e(z), 
$$ 
and  
$$
\langle e^{-{\hbar \over 2} \sum_{j\in\NN} h^{+}[e^{j}] \otimes 
h^{-}[e_{j}] } 
, id 
\otimes q^{((T+U)h^{+})(z)} 
\rangle_{U_{\hbar}\G} =  q^{ - ((T+U)h^{+})(z)},
$$
which follow from direct calculation. 
This finally shows the first part of the lemma. The proof of the second
part is similar. 
\qed\medskip

Let $r\in R$, $\phi\in\cB^{+}$, and let us compute $\langle
F^{-1}[(\Pi\otimes 1)F], id \otimes e[r]\phi\rangle_{U_{\hbar}\G}$. 
We find 
\begin{align} \label{9ab}
\langle
F^{-1}[(\Pi\otimes 1)F], id \otimes e[r]\phi & \rangle_{U_{\hbar}\G}  
= 
\langle
\cR^{-1}[(\Pi\otimes 1)\cR], id \otimes e[r]\phi\rangle_{U_{\hbar}\G} 
\\ \nonumber
& =
\sum \langle \cR^{-1}, id \otimes (e[r]\phi)^{(1)} \rangle_{U_{\hbar}\G} 
\Pi \left( 
\langle \cR, id \otimes (e[r]\phi)^{(2)} \rangle_{U_{\hbar}\G} 
\right)
\\ \nonumber
& = 
\sum S'(e[r]^{(1)}\phi^{(1)})\Pi(e[r]^{(2)}\phi^{(2)}), 
\end{align}
where 
the first equality follows from (\ref{F=R}), and the second one from the
Hopf algebra pairing rules. To show the third one, we remark that
$ ( e[r]\phi )^{(1)}$ belongs to $\cB^{+}$, and apply to it 
Lemma \ref{R=skew}, 1). 

The r.h.s. of (\ref{9ab}) is then 
$\sum S'(e[r]^{(1)}\phi^{(1)})\Pi(e[r]^{(2)}\phi^{(2)})$, but since
$e[r]^{(2)}\in U_{\hbar}\G_{R}$ (see Prop. \ref{additional}), and $\Pi$
is a left $U_{\hbar}\G_{R}$-module morphism, this is
equal to 
$$
\sum S'(e[r]^{(1)}\phi^{(1)})e[r]^{(2)} \Pi(\phi^{(2)})
$$
or 
$$
\sum S'(\phi^{(1)}) S'(e[r]^{(1)}) e[r]^{(2)} \Pi(\phi^{(2)})
$$ 
(because $S'$ is an algebra anti-automorphism of $U_{\hbar}\G$). 
Since $\sum S'(e[r]^{(1)}) e[r]^{(2)}=0$, the r.h.s of (\ref{9ab}) is
equal to zero. By Lemma \ref{prop:ident:pairing}, 1), it then follows that 
$$
F^{-1}[(\Pi\otimes 1)F] \in \cB^{+} \hat\otimes \cB_{R}^{-}, 
$$
that is the first part of (\ref{equiv}). 

The proof of the second part of (\ref{equiv}) is similar: 
let $r\in R$, $\psi\in\cB^{+}$, and let us compute $\langle
[( 1 \otimes \Pi')F] F^{-1}, \psi f[r]\otimes  id \rangle_{U_{\hbar}\G}$. 
The Hopf algebra notation employed now refers to $U_{\hbar}\bar{\G}$.  

We find 
\begin{align} \label{10ab}
  \langle [( 1 \otimes \Pi')F] F^{-1}, \psi & f[r] \otimes id
  \rangle_{U_{\hbar}{\G}} = [( 1 \otimes \Pi')F] F^{-1}, \psi
  f[r]\otimes id \rangle_{U_{\hbar}\bar{\G}} \\ & = \nonumber \langle
  [( 1 \otimes \Pi')\bar \cR^{21}] (\bar \cR^{21})^{-1}, \psi
  f[r]\otimes id \rangle_{U_{\hbar}\bar{\G}} \\ \nonumber & = \sum
  \Pi' \left( \langle \bar \cR^{21}, (\psi f[r])^{(1)} \otimes id
    \rangle_{U_{\hbar}\bar{\G}} \right) \langle (\bar \cR^{21})^{-1},
  (\psi f[r])^{(2)} \otimes id \rangle_{U_{\hbar}\bar{\G}} \\ 
  \nonumber & = \sum \Pi'( (\psi f[r])^{(1)} ) \bar S ( (\psi
  f[r])^{(2)} )
\end{align}
where the first equality follows from Prop. \ref{prop:ident:pairing}, 
the second equality follows from (\ref{F=R}), and third one by the Hopf
algebra pairing rules. To show the last one, we remark that
$ ( \psi f[r])^{(2)}$ belongs to $\cB^{-}$, and apply to it 
Lemma \ref{R=skew}, 2). 

The r.h.s. of (\ref{10ab}) is then 
$\sum \Pi'(\psi^{(1)} f[r]^{(1)} ) 
S (\psi^{(2)} f[r]^{(2)})$, but since
$f[r]^{(1)}\in U_{\hbar}\G_{R}$ (see Lemma \ref{additional}), 
and $\Pi'$
is a right $U_{\hbar}\G_{R}$-module morphism, this is
equal to 
$$
\sum \Pi'(\psi^{(1)}) f[r]^{(1)} 
\bar S (\psi^{(2)} f[r]^{(2)})
$$
or 
$$
\sum \Pi'(\psi^{(1)}) f[r]^{(1)} S( f[r]^{(2)})
\bar S (\psi^{(2)}) 
$$ 
(because $\bar S$ is an algebra anti-automorphism of $U_{\hbar}\G$). 
Since $\sum f[r]^{(1)} S ( \bar f[r]^{(2)} ) =0$, the r.h.s of
(\ref{10ab}) is 
equal to zero. By Lemma \ref{orthBR}, 2), it then follows that 
$$
[( 1 \otimes \Pi')F] F^{-1} \in \cB^{+}_{R} \hat\otimes \cB^{-}, 
$$
that is the second part of (\ref{equiv}). \qed\medskip

\begin{prop} \label{alldecomp}
Any decomposition of $F$ according to (\ref{decomp}) is of the form 
$$
F_{2} = [ (\Pi \otimes 1 ) F ] b , \quad 
F_{1} = a [ (1 \otimes \Pi' ) F ], 
$$
with $a,b\in U_{\hbar}\G_{R}^{\hat\otimes 2}$, such that
$ab=F_{\Pi,\Pi'}$. 
\end{prop}

{\em Proof.} Clear. \qed \medskip

\begin{corollary}
For any left $\cB_{R}^{-}$-module morphism $\bar \Pi$ from $\cB^{-}$ to
$\cB_{R}^{-}$, and any right $\cB_{R}^{+}$-module
morphism $\bar \Pi'$ from $\cB^{+}$ to $\cB_{R}^{+}$, such that $\bar
\Pi(1)=1, \bar \Pi'(1)=1$, we have 
$$
[ (\Pi \otimes 1 ) F ]^{-1}[ (\bar \Pi \otimes 1 ) F ] \in
U_{\hbar}\G_{R}^{\hat\otimes 2}, \quad
[ (1 \otimes \bar \Pi' ) F ][ (1 \otimes \Pi' ) F ]^{-1}\in
U_{\hbar}\G_{R}^{\hat\otimes 2}. 
$$
\end{corollary}

{\em Proof.} This follows from the fact that any $\bar \Pi$, $\bar \Pi'$
yield solutions to (\ref{decomp}), and from the classification of all
such solutions in Prop. \ref{alldecomp}. \qed\medskip

\begin{convention}
The expansion in $\hbar$ of $F$ is $1+\hbar f+o(\hbar)$, in the
notation of Lemma \ref{class-twist}. We may assume that
$\Pi(f[e_{i}])=0$, $\Pi'(e[e_{i}])=0$, for all $i\in \NN$; this implies
that $(1\otimes \Pi')F=1+\hbar f_{1} + o(\hbar)$, and 
$(\Pi\otimes 1)F=1+\hbar f_{2} +
o(\hbar)$, so that $F_{\Pi,\Pi'}=1+o(\hbar)$. In what follows, we will
only consider solutions of (\ref{decomp}), such that 
$F_{1}=1+\hbar f_{1} + o(\hbar)$, and 
$F_{2}=1+\hbar f_{2} + o(\hbar)$; equivalently, the $a$ and $b$ of
Prop. \ref{alldecomp} have the form $1+o(\hbar)$. 
\qed\medskip
\end{convention}

\subsection{Quasi-Hopf structures on $U_{\hbar}\G_{R}$ and $U_{\hbar}\G$}

Let us choose a solution $(F_{1},F_{2})$ of (\ref{decomp}), satisfying
the above requirement. Consider the algebra morphism $\Delta_{R}:
U_{\hbar}\G \to U_{\hbar}\G^{\hat\otimes 2}$, defined as 
\begin{equation} \label{Delta-q-Hopf}
\Delta_{R} = \Ad(F_{1}) \circ \Delta  = \Ad(F_{2}^{-1}) \circ
\bar\Delta;  
\end{equation}
define
\begin{equation} \label{Phi}
\Phi = F_{1}^{23}(1\otimes \Delta)(F_{1}) [F_{1}^{12}(\Delta\otimes
1)(F_{1}) ]^{-1}. 
\end{equation}

\begin{prop} \label{Phi-in-R}
$\Phi$ belongs to $U_{\hbar}\G_{R}^{\hat\otimes 3}$, and even to 
$\cB^+_R \otimes U_\hbar\G_R \otimes \cB^-_R$. 
\end{prop}

{\em Proof.} (\ref{Phi}) makes it clear that $\Phi$ belongs to
$U_{\hbar}\G^{\hat\otimes 2}\hat\otimes U_{\hbar}\G_{R}$. It can be
rewritten as 
$$
\Phi = (F_{2}^{-1})^{23}(1 \otimes \bar\Delta)(F_{2}^{-1})
[ (F_{2}^{-1})^{12}(\bar\Delta \otimes 1)(F_{2}^{-1}) ]^{-1}; 
$$
it follows that $\Phi\in U_{\hbar}\G_{R}\hat\otimes 
U_{\hbar}\G^{\hat\otimes 2}$. Finally, $\Phi$ can also be written as 
$$
\Phi = [(1\otimes \Delta_{R})F_{1}](F_{2}^{-1})^{23}F^{23}
(\Delta \otimes 1)(F^{-1})(\Delta\otimes 1)(F_{2})
(F_{1}^{-1})^{12}.
$$
But $(\Delta\otimes 1)(F^{-1})$ belongs to $U_{\hbar}\N_+ \otimes 
U_{\hbar}\G_+ \otimes U_{\hbar}\N_-$, by Prop. \ref{additional}. 
Therefore, $\Phi$ belongs to $U_{\hbar}\G \otimes (U_{\hbar}\G_R
U_{\hbar}\G_+ U_{\hbar}\G_R) \otimes U_{\hbar}\G$.

$\Phi$ can again be rewritten as 
$$ 
\Phi = F_2^{-1(23)} (1\otimes \bar\Delta)(F_1) (1\otimes
\bar\Delta)(F^{-1}) F^{(12)} (\Delta\otimes 1)(F_2) F_1^{-1(12)}. 
$$ 
Prop. \ref{additional} now implies that $(1\otimes \bar\Delta)(F_1)$
belongs to $U_{\hbar}\G \otimes U_{\hbar}\G_- \otimes U_{\hbar}\G$.
Therefore, $\Phi$ belongs to $U_{\hbar}\G \otimes (U_{\hbar}\G_R
U_{\hbar}\G_- U_{\hbar}\G_R) \otimes U_{\hbar}\G$.

By the PBW results of sect. \ref{PBW}, the intersection of
$U_{\hbar}\G_R U_{\hbar}\G_+ U_{\hbar}\G_R$ and $U_{\hbar}\G_R
U_{\hbar}\G_- U_{\hbar}\G_R$ is reduced to $U_{\hbar}\G_R$.
Therefore, $\Phi$ belongs to $U_{\hbar}\G \otimes U_{\hbar}\G_R
\otimes U_{\hbar}\G$.  \qed\medskip

Let 
\begin{equation} \label{uR}
u_{R} = m(1\otimes
S)(F_{1}),
\end{equation} 
with  $m$ the multiplication of $U_{\hbar}\G$. 

\begin{theorem}
The algebra $U_{\hbar}\G$, endowed with the coproduct $\Delta_{R}$,
associator $\Phi$, counit $\varepsilon$, antipode $S_{R} = \Ad(u_{R})
\circ S$, respectively
defined in 
(\ref{Delta-q-Hopf}), (\ref{Phi}), (\ref{counit}), 
(\ref{S-efhKD}), (\ref{S-ef}), (\ref{uR}), and $R$-matrix
\begin{equation} \label{R:mx:fin}
\cR_{R} =  [ a^{21}(\Pi'\otimes 1)(F^{21}) ] 
q^{D\otimes K} q^{{1\over 2}\sum_{i\in \NN}h^{+}[e^{i}]\otimes
h^{-}[e_{i}]}
[(\Pi \otimes  1)(F)F_{\Pi,\Pi'}a^{-1}], 
\end{equation}
is a quasi-triangular quasi-Hopf algebra. $U_{\hbar}\G_{R}$ is a
sub-quasi-Hopf algebra of it. Moreover, $\cR_R$ belongs to $U_{\hbar}\G_R
\hat\otimes U_{\hbar}\G$. 
\end{theorem}

{\em Proof.} The statement on $U_{\hbar}\G$ follows directly from
\cite{D}, sect. 1, rem. 5. 
That $U_{\hbar}\G_{R}$ is a sub-quasi-bialgebra of $U_{\hbar}\G$ 
follows from Prop. \ref{why-decomp} and Prop. \ref{Phi-in-R}. Let us now
show that $S_{R}$ preserves $U_{\hbar}\G_{R}$. 

We have $\sum_{i}x_{i}S_{R}(x'_{i}) = \varepsilon(x)$, for $x\in
U_{\hbar}\G$, where $\Delta_{R}(x) = \sum_{i} x_{i} \otimes x'_{i}$. Let
$(m_{\al})$ be a basis of $U_{\hbar}\G$ as a left
$U_{\hbar}\G_{R}$-module, with $m_{0} =1$. Set $S_{R} =
\sum_{\al}S_{R}^{(\al)} m_{\al}$, with $S_{R}^{(\al)}$ some linear map
from $U_{\hbar}\G_{R}$ to itself. Then for $\al\neq 0$, 
$\sum_{i}x_{i}S_{R}^{(\al)}(x'_{i}) =0$. Let us show that this implies
that $S_{R}^{(\al)} =0$. 

Assume that for some $\al$, $S_{R}^{(\al)}$ is not $0$. Dividing it by the
largest possible power of $\hbar$, we may assume that its classical
limit $S_{R,cl}^{(\al)}$ is non-zero. $S_{R,cl}^{(\al)}$ is then a map
from $U\G_{R}$ to itself, such that $\sum_{i} y_{i} S_{R,cl}^{(\al)}
(y'_{i} ) =0$ for any $y\in U\G_{R}$, where $\Delta(y) = \sum_{i}
y_{i}\otimes y'_{i}$. We then check by induction on the degree of $y$
that $S_{R,cl}^{(\al)} =0$, a contradiction. 

So $S_{R} = S_{R}^{(0)}$, and $S_{R}$ preserves $U_{\hbar}\G_{R}$. 

That $\cR_R$ belongs to $U_{\hbar}\G_R \hat \otimes U_{\hbar}\G$
follows clearly from (\ref{R:mx:fin}).  \hfill \qed

\subsection{Adelic algebras} \label{adal}

In \cite{D}, Drinfeld also defined an adelic version of the Manin pair
$(\A\otimes k, \A \otimes R)$. Let $\AAA$ be the ring of adeles of $X$
and $\CC(X)$ be the field of meromorphic functions of $X$. Let us define
on $\AAA$ the scalar product $\langle f,g \rangle_{\AAA} = \sum_{x\in
X}\res_{x}(fg\omega)$. Endow $\A \otimes \AAA$ with the scalar product
of the Killing form of $\A$ with $\langle, \rangle_{\AAA}$. The Lie
algebra $\A \otimes \CC(X)$ is then a Lagrangian subspace of $\A \otimes
\AAA$; the pair $(\A\otimes \AAA, \A \otimes \CC(X))$ then forms a Manin
pair. 

It is easy to double extend it as in section \ref{mp}.  The
construction of quasi-Hopf algebras presented here then can be applied
to yield an ``adelic'' quasi-Hopf algebra quantizing this Manin pair.

\subsection{Quantum Weyl group action} \label{qwga}

The Weyl group $W$ of $\A$
naturally maps to the group of automorphisms
of the Manin pair $(\G,\G_{R})$. There is an algebra automorphism
$w$ of $U_{\hbar}\G$, deforming the action of the nontrivial element of
$W$. It is defined by the rules 
$$
( w\cdot e)(z) = -\left( q^{K\pa}(q^{h^{-}}f)\right)(z), \quad 
(w\cdot f) (z) = -  e(z) q^{-((T+U)h^{+})(z)},
$$
$$
w(h^{+}[r])= - h^{+}[r], \quad w(h^{-}[\la])= - h^{-}[\la], \quad
w(D)=D, w(K)= K, 
$$
where $r\in R, \la\in \La$ and $(w\cdot x)(z) = \sum_{i\in \ZZ}
w(x[\eps^{i}])\eps_{i}(z)$, $x=e,f$.  

Note that $w$ does not preserve $U_{\hbar}\G_{R}$, and $w^{2}\neq 1$.
\medskip

\section{Analogues and generalizations of $U_{\hbar}\G$} \label{agen}

\subsection{Replacing $q^{K\pa}$ by a general automorphism}
\label{galaut} 

The algebras $U_{\hbar}\G$ presented in section \ref{ug} admit the
following generalizations. 
Let $\sigma$ be a ring automorphism of $k$, commuting with $\pa$
and preserving $R$ and $\langle, \rangle_{k}$. Then we can
form an algebra $U_{\hbar}\G_{k,\sigma}$ with the same generators as
$U_{\hbar}\G$ (except $K$), replacing in all relations $q^{K\pa}$ by
$\sigma$. For example, (\ref{h-h}) becomes 
$$
[h^{+}[r], h^{-}[\la]] = {2\over \hbar}  \langle (1-\sigma^{-1})r, \la
\rangle_{k}, 
$$
etc. Expressing the action of $\Ad(q^{KD})$ on $U_{\hbar}\G$ and
replacing in the resulting formulas, $q^{K\pa}$ by $\sigma$, we obtain
an (outer) automorphism $\Sigma$ of $U_{\hbar}\G_{k,\sigma}$. The
formulas for $\Sigma$ are 
$$
\Sigma(q^{ - ((T+U)h^{+})(z)} x(z)) = \sigma^{-1}(q^{ - ((T+U)h^{+})(z)}
x(z)), \quad \Sigma(D)= D, 
$$
$$
\Sigma(h^{+}[r])=h^{+}[\sigma(r)], \quad
\Sigma(h^{-}[\la])=h^{-}[(\sigma(\la))_{\La}] + h^{+}[\sigma((T+U)\la) -
(T+U)((\sigma(\la))_{\La})],  
$$
$x=e,f$; in the r.h.s. of the first formula, $\sigma^{-1}$ is applied to
the function part. 

Note that $U_{\hbar}\G_{R}$ is a subalgebra of $U_{\hbar}\G_{k,\sigma}$,
and that $\Sigma$ restricts to an automorphism of $U_{\hbar}\G_{R}$. 

If $\sigma$ is finite, then we can find some $\tau$ such that the
formulas defining the action of $\Sigma$ on $x(z)$ are simply
$\Sigma(x(z)) = \sigma^{-1}(x(z))$, $x=e,f$. Indeed, in that case the
equation 
$$
(\sigma\otimes \sigma) \tau - \tau + \sum_{i\in \NN}
[T(\sigma((\sigma^{-1}(e_{i}))_{\La})) - T(e_{i}) ] \otimes e^{i} =0
$$
can easily be solved. 

If moreover $\sigma(\La) = \La$, then $\Sigma$ coincides with
$\sigma^{-1}$ also on the generating series $h^{\pm}(z)$. 

\begin{remark} 
As in Remark \ref{frob}, the algebra $U_{\hbar}\G_{k,\sigma}$ (without
$D$) can be generalized with $k$, $\langle 1, \cdot\rangle_{k}$, 
$R\subset k$ and $\sigma$ 
replaced by an arbitrary Frobenius algebra $(k_{0}, \theta)$ 
with a Lagrangian subalgebra
and an automorphism, preserving the scalar product and the subalgebra; 
and the full algebra $U_{\hbar}\G_{k,\sigma}$ can be generalized to the
case where $k_{0}$ is endowed in addition with a derivation $\pa_{0}$
commuting with the automorphism, and such that $\theta\circ \pa_{0} =0$. 
\qed
\end{remark}

\begin{remark} It seems difficult to combine the above action with the
quantum Weyl group action of section \ref{qwga} to give
quantizations of more general ``twisted'' Manin pairs of
Drinfeld. In that situation, $X$ is endowed with an involution $\sigma$
preserving $\omega$ and $S$, and the Manin pair is defined in the
algebra $({\mathfrak{sl}}_{2}\otimes k)^{\ZZ/2\ZZ}$, where the action of
$\ZZ/2\ZZ$ is by the tensor product of the Weyl group
action with $\sigma$. The difficulty is that $\Sigma$ has finite
order whereas $w$ has infinite order. If after conjugation, $w$ could be
brought to a ``correct'' form (that is, preserving $U_{\hbar}\G_{R}$), it
might happen that the procedure described here applies. 
\end{remark}

\subsection{Discrete analogues} \label{disc}

Let us formally express the algebra relations of $U_{\hbar}\G$, using
the new generating series $K^{+}(z) = q^{((T+U)h^{+})(z)}$ and
$K^{-}(z)=q^{h^{-}(z)}$. Let us replace the exression $q^{2\sum_{i\in
\NN}((T+U)e_{i})(z)e^{i}(w)}$ by $q(z,w)$. Using (\ref{id-tau}), we
formally derive the equation 
\begin{equation} \label{antisymm-q}
q(z,w) q(w,z) =1; 
\end{equation}
let us also note $(q^{K\pa}f)(z)$ as $f(\sigma^{-1}(z))$. Then the
formulas presenting $U_{\hbar}\G$ become
\begin{equation} \label{K-K}
(K^{+}(z),K^{+}(w))=1, \quad
(K^{+}(z),K^{-}(w))={{q(z,w)}\over{q(z,\sigma(w))}},
\end{equation}
\begin{equation} \label{K-K-}
(K^{-}(z),K^{-}(w))={{q(\sigma(z),\sigma(w))}\over{q(z,w)}},
\end{equation}
\begin{equation}\label{K-e}
(K^{+}(z),e(w))= q(z,w), \quad (K^{-}(z),e(w))= q(w, \sigma(z)), 
\end{equation}
\begin{equation}\label{K-f}
(K^{+}(z),f(w))= q(w,z), \quad (K^{-}(z),f(w))= q(z,w), 
\end{equation}
\begin{equation} \label{dsc-vx}
(e(z), e(w)) = q(z,w), \quad (f(z), f(w)) = q(w,z), 
\end{equation}
\begin{equation} \label{dsc-e-f}
[e(z),f(w)] = \delta_{z,w}K^{+}(z) -\delta_{z, \sigma(w)}K^{-}(w)^{-1}; 
\end{equation}
we used the standard notation $(a,b)$ for the group commutator
$aba^{-1}b^{-1}$; we also multiplied $f(z)$ by $\hbar$, and replaced
$\delta(z,w)$ by $\delta_{z,w}$. 
We have also the trivial relations 
\begin{equation} \label{triv}
K^{\pm}(z)K^{\pm}(z)^{-1}=K^{\pm}(z)^{-1}K^{\pm}(z)=1.
\end{equation}

Generators $K^{\pm}(z), K^{\pm}(z)^{-1}, e(z)$ and $f(z)$
and relations (\ref{K-K}), (\ref{K-K-}), (\ref{K-e}), (\ref{K-f}),
(\ref{dsc-vx}), (\ref{dsc-e-f}) can be thought of as presenting a
complex algebra
$A(Z,\sigma, q)$ defined from the data of a discrete set $Z$, a map
$\sigma:Z\to Z$, and a function $q:Z^{2}\to \CC^{\times}$, satisfying
(\ref{antisymm-q}). 
It is easy to see that a basis for $A(Z,\sigma,q)$ is formed by the
family $\prod_{z\in Z} e(z)^{\eps_{z}} \prod_{z\in Z}
K^{+}(z)^{\kappa^{+}_{z}} \prod_{z\in Z} K^{-}(z)^{\kappa^{-}_{z}}
\prod_{z\in Z} f(z)^{\eta_{z}}$, $\eps_{z},\eta_{z}\in \NN$,
$\kappa^{\pm}_{z}\in \ZZ$. 

The quantum Weyl group action of section \ref{qwga} then has the
following discrete analogue. Assume $\sigma$ to be invertible. Then
there is an automorphism $w_{Z}$ of 
$A(Z,\sigma,q)$, defined by the formulas 
$$
w_{Z}(e(z)) = -K^{-}(\sigma^{-1}(z)) f(\sigma^{-1}(z)), \quad
w_{Z}(f(z)) = - e(z) K^{+}(z)^{-1}, 
$$
$$
w_{Z}(K^{\pm}(z)) = K^{\pm}(z)^{-1}. 
$$

In the case where $\sigma= id_{Z}$, the discrete analogue of the
coalgebra structure of $U_{\hbar}\G$ is then given by 
$$
\Delta(K^{\pm}(z))= K^{\pm}(z)\otimes K^{\pm}(z), \quad 
\Delta(e(z)) = e(z)
\otimes K^{+}(z) + 1 \otimes e(z),
$$
$$
\Delta(f(z))= f(z)\otimes 1 +
K^{-}(z)^{-1} \otimes f(z).
$$
The subalgebras $A_{+}(Z, id_Z, q)$ and $A_{-}(Z, id_Z,q)$ of
$A(Z,id_Z,q)$ generated by the $e(z)$ and $K^{+}(z)^{\pm1}$, resp. by
the $f(z)$ and $K^{-}(z)^{\pm 1}$ then form Hopf subalgebras of
$A(Z,id_{Z},q)$. If $A_{+}(Z, id_{Z},q)$ is given the opposite
coproduct, they are in duality, the pairing being defined by 
$\langle e(z),f(w) \rangle =\delta_{z,w}$, 
$\langle K^{+}(z)^{\eps},K^{-}(w)^{\eps'} \rangle =q(z,w)^{\eps\eps'}$,
$\eps,\eps'=1$ or $-1$. 

\begin{remark}
It is also natural to ask whether the formal series in $\hbar,z$ and $w$
given by 
$\exp(2\hbar\sum_{i}((T+U)h^{+})(z)e_{i}(w))$ 
has an analytic prolongation. It
could then happen that the relations defining $A(Z,\sigma,q)$ can be
represented over $\CC$
in a weak sense -- as relations between analytic prolongations of 
matrix coefficients of operators acting on highest weight modules. 
\medskip
\end{remark}


\begin{remark} \label{correction}
  In the situation where $k$ and $R$ are replaced by a finite
  dimensional Frobenius algebra and a maximal isotropic subring of it,
  an expression for $F$ is \begin{equation}
    \label{F=exp} F = \exp(\hbar \sum_i e[\eps^i] \otimes
    f[\eps_i]),\end{equation} for $(\eps^i),(\eps_i)$ two dual bases
  of $k$. As it was pointed out in \cite{DK}, this result is no longer
  true when $k$ is infinite-dimensional. For example, such results as
  the commutativity of $\sum_i e[\eps^i]\otimes f[\eps_i]$ with $e(z)
  \otimes K^+(z) + 1 \otimes e(z)$ cease to be true in the case where
  $(k,R)$ are associated to a curve with marked points.  This is
  because, by vanishing properties of $1 + a_0 \psi_+$ (see
  \cite{Some}), the defining relations only imply an identity
  $$(z-q^{-\pa}w) [(e\otimes f)(z), (e\otimes q^{(T+U)h^+})(w)]= 0,$$
  so that $$[(e\otimes f)(z), (e\otimes q^{(T+U)h^+})(w)]= A(z)
  \delta(z,q^{-\pa}w),$$ for some field $A(z)$, so that $[\sum_i
  e[\eps^i] \otimes f[\eps_i], (e\otimes q^{(T+U)h^+})(w)] =
  A(q^{-\pa} w)$.

  In \cite{DK}, an expression for $F$, well-defined up to order
  $2$, was proposed in the quantum affine case and checked up to
  that order.

  It would be interesting to obtain expressions for $F$ is the
  framework of \cite{FO}; one could expect to check their coincidence
  with those of \cite{KT}.

  An earlier version of the present work used a (wrong) generalization
  of (\ref{F=exp}) to the infinite-dimensional setting; after that,
  the works \cite{EF1,EF2} were completed, relying on this work.
  However, after $F$ is defined by Def. \ref{def:def:F}, all results
  of that version are correct, except those involving commutation
  relations of $\sum_i e[\eps^i]\otimes f[\eps_i]$.  This shows that
  the only corrections to \cite{EF1,EF2} are to replace the definition
  of $F$ based on the generalization of (\ref{F=exp}) by
  Def. \ref{def:def:F}.
\end{remark}

\section*{Appendix: maximal isotropy of rings}

Let $\AAA$ be the adeles ring of
$X$, and $\CC(X)$ be the field of functions over $X$. Define on $\AAA$
the bilinear pairing $\langle, \rangle_{\AAA}$ by 
$$
\langle f, g\rangle_{\AAA}=\sum_{x\in X}\res_{x}(fg\omega).
$$
In \cite{D}, Drinfeld made the following statement. 

\begin{lemma} \label{statement-K}
$\CC(X)$ is maximal isotropic in $\AAA$ w.r.t. $\langle, \rangle_{\AAA}$. 
\end{lemma}

Below we prove this and the similar statement

\begin{lemma} \label{statement-R}
$R$ is maximal isotropic in $k$ w.r.t. $\langle, \rangle_{k}$. 
\end{lemma}

Recall first
the duality theorem (\cite{Serre}, II-8, thm. 2). 
Let $D$ be any divisor
on $X$, and $\Omega(D)$ be the space of all meromorphic forms $\omega$
equal to zero or such that their divisor is $\ge D$. Let on the other
hand, $\AAA_{\ge -D}$ be the space of adeles with divisor $\ge -D$. Then
$\langle, \rangle_{\AAA}$ induces a non-degenerate pairing 
$$
\Omega(D)\times \left(\AAA/(\AAA_{\ge -D}+\CC(X))\right) \to \CC. 
$$

\noindent{\em Proof of Lemma \ref{statement-K}.}
The isotropy of $\CC(X)$ follows from the residue formula. 
Let $\Omega$ be the space of all meromorphic forms on $X$, and let us 
now show that the pairing 
\begin{equation} \label{adelic-pairing}
\Omega \times (\AAA/\CC(X)) \to \CC
\end{equation}
is also non-degenerate. Let $f\in \AAA/\CC(X)$ have vanishing pairing with
$\Omega$. Then for any divisor $D$, the pairing of 
its image in $\AAA/(\AAA_{\ge -D}+\CC(X))$ with any element of $\Omega(D)$ is
zero, which means that $f$ belongs to $\AAA_{\ge -D}/(\AAA_{\ge -D}\cap
\CC(X))$ for any $D$, and so is zero. 

The lemma now follows from the non-degeneracy of (\ref{adelic-pairing}). 
\qed
\medskip

\noindent
{\em Proof of Lemma \ref{statement-R}.}
Let for any divisor $\bar D$ with support in $S$, $k_{\ge \bar{D}}$
be the space of elements of $k$ with divisor $\ge \bar{D}$. 

\begin{lemma} \label{A-k}
Let $D_{0}$ be a divisor of $X$, supported in $S$. Then 
the mappings $\AAA \to k$ induces an isomorphism of 
$\AAA/(\AAA_{\ge -D_0}+\CC(X))$ with $k/ (k_{\ge -D_0}+R)$. 
\end{lemma}

{\em Proof.}
Let $D=n(\sum_{s\in S}s)
$. For $n$ large enough, $D \ge ( \omega_0)$ 
and $\Omega(D)=0$; by the duality theorem, this implies that 
$ \AAA/( \AAA_{ \ge -D}+\CC(X))
=0$, and so $ \AAA= \AAA_{\ge -D}+\CC(X)$. 
Let $D_0$ be a divisor 
$\le D$, then $\AAA_{\ge -D_0}\subset \AAA_{\ge -D}$, so 
\begin{align*}
\AAA/( \AAA_{\ge -D_0}+\CC(X))
& =
( \AAA_{ \ge -D}+\CC(X))/( \AAA_{ \ge -D_0}+\CC(X))
\\ & =
\AAA_{ \ge -D}/( \AAA_{ \ge -D_0}+(\CC(X) \cap  \AAA_{\ge -D})). 
\end{align*}
Let for any $n$, 
$Q_n=\AAA_{ \ge -n(\sum_{s\in S}s)}/\left( \AAA_{ \ge -D_0}
+(\CC(X) \cap  \AAA_{\ge -n(\sum_{s\in S}s)})\right)$. 
For the same $n$
as above, the natural maps $Q_n \to Q_{n+1}$ 
are isomorphisms. On the other hand, we have a map 
$Q_n \to k/(k_{\ge -D_0}+R)$; 
its kernel is the 
set of adeles $\ge -n(\sum_{s\in S}s)$, $\ge -D_0$, and in $R$, 
so is zero. So 
$Q_n \to k/(k_{\ge -D_0}+R)$ is injective; it is also 
surjective, as we can see by composing it with a suitable 
$Q_n \to Q_m$, $m$ large enough, so it is an isomorphism. 
\qed\medskip

>From the duality theorem now follows that for $D_{0}$ like in Lemma 
\ref{A-k}, the residue pairing 
$$
\Omega(D_0) \times \left( k/ (k_{\ge -D_0}+R) \right) \to \CC
$$
is non-degenerate.
Let us specialize $D_{0}$ to 
$D_{0,m}=-m(\sum_{s\in S}s)$, then the natural maps induce 
an inductive system $\Omega_{0,m}\subset \Omega_{0,m'}$ and a projective
system $k/
(k_{\ge -D_{0,m'}}+R) \to k/
(k_{\ge -D_{0,m}}+R) $, $m\le m'$, compatible with
the duality. It follows that the induced pairing between the inductive
and projective limits is non-degenerate. Since 
$\cup_{m\ge 0}\Omega(D_{0,m})$ is the set of forms on $X$ regular outside
$S$, and $\omega$ has neither zeros nor poles outside $S$, this space
is equal to $\omega R$. On the other hand, 
$\limm_{m\ge 0}k/
(k_{\ge -D_{0,m}}+R) =
k/R $. We conclude that the pairing 
$$
\omega R \times (k/R) \to \CC
$$
induced by the residue is non-degenerate, whence the lemma. 
\qed\medskip

\begin{remark} Another proof of Lemma \ref{statement-R} can be obtained
following \cite{D}, sect. 2. \end{remark}

\end{document}